\documentclass[reqno]{article}
\usepackage{amsmath,amstext,amssymb,amscd,bbm,amsthm}
\usepackage{hyperref}

\numberwithin{equation}{section} \theoremstyle{plain}
\newtheorem{lemma}{Lemma}[section]
\newtheorem{theorem}[lemma]{Theorem}
\newtheorem{corollary}[lemma]{Corollary}
\newtheorem{proposition}[lemma]{Proposition}
\newtheorem{definition}[lemma]{Definition}
\newtheorem{remark}{Remark}[section]
\newtheorem{example}{Example}[section]

\def\a{\alpha}
\def\b{\beta}

\def\g{\gamma}

\def\om{\omega}

\def\G{\Gamma}

\def\m{\mu}

\newcommand{\innerprod}[1]{\left\langle#1\right\rangle}
\newcommand{\norm}[1]{\left\|#1\right\|}
\newcommand{\abs}[1]{\left|#1\right|}
\newcommand{\floor}[1]{\left\lfloor#1\right\rfloor}  
\newcommand{\brac}[1]{\left\{#1\right\}}

\newcommand{\R}{\mathbb R} 
\newcommand{\Z}{\mathbb Z} 
\newcommand{\E}{\mathbb E}  
\newcommand{\N}{\mathbb N}  
\newcommand{\Tor}{\mathbb T}

\newcommand{\sC}{\mathcal C}

\renewcommand{\Tor}{\mathbb T}

\newcommand{\X}{\mathcal X}
\newcommand{\Y}{\mathcal Y}
\newcommand{\I}{\mathcal I}

\newcommand{\sZ}{\mathcal Z} 
\newcommand{\sG}{\mathcal G} 
 
\newcommand{\bfs}{\textbf{s}}  
\newcommand{\bfu}{\textbf{u}}
\newcommand{\bfv}{\textbf{v}}

\newcommand{\bfn}{\textbf{n}}
\newcommand{\bfx}{\textbf{x}}

\newcommand{\bfRN}{\textbf{R}_N}
\newcommand{\folN}{\brac{\Phi_N}_{N\in \N}}

\newcommand{\limsupN}{\limsup_{N \rightarrow \infty}}
\newcommand{\limN}{\lim_{N \rightarrow \infty}}
\newcommand{\avgNfol}{\frac{1}{m(\Phi_N)} \int_{\Phi_N}}
\newcommand{\avgs}[1]{\frac{1}{#1} \sum_{n=0}^{#1-1}}
\newcommand{\avgi}[1]{\frac{1}{#1} \int_{0}^{#1}}

\newcommand{\avgNsub}[1]{\frac{1}{N_{#1}} \sum_{n_{#1}=0}^{N_{#1}-1}}
\newcommand{\avgRd}{\frac{1}{R_1}\int_{0}^{R_1} \ldots \frac{1}{R_d}\int_{0}^{R_d}}

\newcommand{\bm}[1]{\mbox{\boldmath{$#1$}}}
\newcommand{\seminorm}[1]{ \left\vert\!\left\vert\!\left\vert #1\right\vert\!\right\vert\!\right\vert}

\begin{document}

\title{Multiple ergodic averages for flows and an application.}
\author{Amanda Potts\footnote{Department of Mathematics, Northwestern University, 2033 Sheridan Road, Evanston, IL 60208-2730, USA, apotts@math.northwestern.edu}}

\date{\today}

\maketitle

\begin{abstract}  We show the $L^2$-convergence of continuous time ergodic averages of a product of functions evaluated at return times along polynomials.  These averages are the continuous time version of the averages appearing in Furstenberg's proof of Szemer\'edi's Theorem.  For each average we show that it is sufficient to prove convergence on special factors, the Host-Kra factors, which have the structure of a nilmanifold.  We also give a description of the limit.  In particular, if the polynomials are independent over the real numbers then the limit is the product of the integrals.  We further show that if the collection of polynomials has ``low complexity'', then for every set $E$ of real numbers with positive density and for every $\delta >0$, the set of polynomial return times for the ``$\delta$-thickened'' set $E_{\delta}$ has bounded gaps.  We give bounds for the flow average complexity and show that in some cases the flow average complexity is strictly less than the discrete average complexity.
\end{abstract}

\section{Introduction.}\label{bsec1}\subsection{Multiple convergence for flows.}\label{bsec1.1}  Furstenberg's groundbreaking proof of Szemer\'edi's theorem via ergodic theory gave rise to many interesting avenues of research.  Of particular importance, it established the connection between recurrence properties of subsets of $\N$ and the limiting behavior of certain associated multiple ergodic averages.  In this paper we focus on the natural analogues of some of these results for multiple ergodic averages along flows.  Let $m$ denote Lebesgue measure on $\R^d$, $d\in\N$.  We show:
\begin{theorem} \label{mainconvg}  Let $\{T_t\}_{t\in\R}$ be a measure preserving flow on a Lebesgue space $(X,\X,\m)$ and let $\{p_1,\ldots, p_k\colon \R^d  \rightarrow \R\}$ be any collection of polynomials.  Then for any $k \in \N$, and $f_1, \ldots, f_k \in L^\infty(\m)$, 
\begin{equation}\label{maineq}  \frac{1}{R_1}\int_{0}^{R_1} \ldots \frac{1}{R_d}\int_{0}^{R_d} f_1 \circ T_{p_1(\bfs)} \cdot \ldots \cdot  f_k \circ T_{p_k(\bfs)}\ d\bfs\end{equation}
converges in $L^2(\m)$ as $R_1,\ldots,R_d \rightarrow \infty$.
\end{theorem}
It is known in the discrete case that for polynomials $\Z^d \rightarrow \Z$, the multiple polynomial averages for a single ergodic transformation converge in $L^2(\m)$, with results given in \cite{Furst1,CL,FWa,HK,HK2,L1}.

In this paper we also describe the limit of \eqref{maineq}.  If $\brac{p_1,p_2,\ldots,p_k}$ is a family of polynomials which are independent over the real numbers, we show that the average \eqref{maineq} converges to the product of the integrals:
\begin{theorem}\label{btheorem1.2}  Suppose $\{T_t\}_{t\in\R}$ is an ergodic measure preserving flow on a Lebesgue space $(X,\X,\m)$, the family of polynomials $\brac{p_1,\ldots,p_k\colon \R^d \rightarrow \R}$ is $\R$-independent, and $f_1,\ldots,f_k\in L^{\infty}(\m)$.  Then as $R_1,\ldots,R_d \rightarrow \infty$,
\begin{equation*} \avgRd f_1 \circ T_{p_1(\bfs)}\cdot  f_2\circ T_{p_2(\bfs)} \cdot \ldots \cdot  f_k \circ T_{p_k(\bfs)} \,d\bfs \,\end{equation*}
 converges in $L^2(\m)$ to 
 \begin{equation*}\int f_1 \,d\m \cdot \int f_2 \,d\m\cdot \ldots \cdot \int f_k \,d\m.\end{equation*} \end{theorem}The discrete version of Theorem \ref{btheorem1.2} was proved in \cite{FK2}.

We also give a formula for the limit of \eqref{maineq} when $p_1,\ldots,p_k$ are not necessarily independent (see discussion in Section \ref{bsec5.2}).  In the discrete setting, an explicit formulation of the limit is given for various cases in \cite{Z2,Fran,Le4}.  In the setting of a flow, the extra level of connectedness in the underlying space allows us to give an explicit description of the limit in general.  

\subsection{Optimal lower bounds.}\label{bsec1.2}

Suppose $f_1=\ldots =f_k = \bm{1}_A$ for some measurable set $A$.  In this situation, Theorem \ref{btheorem1.2} shows that the best lower bound we could expect for \eqref{maineq} is $\m(A)^k$.  We know that in general the limit is not $\m(A)^k$ (see Section \ref{bsec5.2} for a counter-example; see \cite{Z2,Le4,BHK,Fran} for counterexamples in the discrete case).  However, we show that under certain conditions the average is frequently greater than $\m(A)^k-\varepsilon$ for every $\varepsilon >0$.  We say a set $S \subseteq \R^d$ is \emph{\textbf{syndetic}} if there exists a compact set $C \subset \R^d$ such that $\R^d=C + S$.  We show that for collections of polynomials with \emph{complexity} 0 or 1 (see Section \ref{bsec5.3} for the definition), the optimal lower bound is reached for a syndetic set of times:
\begin{theorem}\label{btheorem1.3}  Suppose $\{T_t\}_{t\in\R}$ is an ergodic measure preserving flow on a Lebesgue space $(X,\X,\m)$, $A\in \X$ with $\m(A)>0$, and $\{ p_1,\ldots,p_k  \colon \R^d \rightarrow \R\}$ are polynomials with $p_i(0)=0$ for $i=1,\ldots,k$.  If $\{ p_1,\ldots,p_k  \}$ has complexity 0 or 1 then for every $\varepsilon >0$ the set 
$$ \{ \bfs \in \R^d \colon \m(A \cap T_{p_1(\bfs)}A \cap\ldots \cap T_{p_k(\bfs)}A) \geq \m(A)^{k+1}-\varepsilon   \}  $$
is syndetic. 
\end{theorem}
We note that a family of polynomials has complexity 0 if and only if it is $\R$-independent.  Some examples of families with complexity 1 are $\{t,t^2,t+t^2\}$, $\{t,2t,t^2\}$, and $\{t,t^2,t^3,t+t^2+t^3\}$.

The discrete time version of Theorem \ref{btheorem1.3} for polynomials of the form $\{n,2n\}$ and $\{n,2n,3n\}$ was given by Bergelson, Host, and Kra in \cite{BHK}, and was generalized by Frantzikinakis in \cite{Fran} to include all collections of three polynomials of \emph{Weyl complexity} 1 or 2 (see \cite{BLL} for the definition).  For the discrete case, it is known that the optimal lower bound is \emph{not} reached for the polynomial family $\{n,2n,3n,4n\}$ (see \cite{BHK}).  We note that here the discrete and continuous versions differ, as there exist collections of three polynomials which have Weyl complexity 3, but have complexity 1.  One such collection is $\{n,2n,n^2\}$, for which the discrete version of Theorem \ref{btheorem1.3} is likely to fail \cite{Fran} (this is currently unknown).  

We give a family of polynomials with complexity 2 which achieves the optimal lower bound:
\begin{theorem}\label{btheorem1.4}  Suppose $\{T_t\}_{t\in\R}$ is an ergodic measure preserving flow on a Lebesgue space $(X,\X,\m)$, $A\in\X$ with $\m(A)>0$, and $l,m\in\N$.  If $p \colon \R^d \rightarrow \R$ is a polynomial with $p(0)=0$, then for every $\varepsilon >0$, the set 
\begin{equation*} \{ \bfs \in \R^d \colon \m(A \cap T_{l p(\bfs)}A \cap T_{m p(\bfs)}A \cap T_{(l+m) p(\bfs)} A) \geq \m(A)^4-\varepsilon   \}   \end{equation*}
is syndetic.
\end{theorem}
It is unknown whether Theorem \ref{btheorem1.3} holds for families of complexity 2.  In the discrete case, $\{2n,3n,4n\}$ is a family of complexity 2 for which the discrete version of Theorem \ref{btheorem1.3} is likely to fail \cite{Fran}.  The discrete analog of Theorem \ref{btheorem1.4} was given in \cite{Fran}.

\subsection{Application.}\label{bsec1.3}  Just as Furstenberg used ergodic results \cite{Furst1} to derive Szemer\'edi's Theorem, we are able to derive combinatorial results from our study of continuous time averages.  In particular, given a sufficiently large subset $E \subseteq \R$, we ask which types of configurations are guaranteed to lie arbitrarily close to $E$.  Let us make this question more precise.  The \emph{\textbf{upper Banach density}} of a subset $E \subseteq \R$ is the quantity $$D^*(E)= \limsup_{(N-M) \rightarrow \infty} \frac{m\bigl(E \cap [M,N]\bigr)}{(N-M)}.$$  For $\delta>0$, we write $E_{\delta}\colon = \{v \in \R \colon \text{dist}(v,E)<\delta \}=\{v \in \R \colon \abs{v-e}<\delta\, \,\text{for some}\,e\in E \}$.  If $E \subseteq \R$ with $D^*(E)>0$, we are interested in paths $\{a_1(t),\ldots,a_k(t)\}_{t\in \R}$ which have the property that for each $\delta >0$ there exists $x,t_0\in\R$ with $x+a_1(t_0),\ldots,x+a_k(t_0)\in E_{\delta}$.  For example, it is shown in \cite{Z} that given $\{\a_1,\ldots,\a_k\}\subset\R$ and $\delta>0$, there exists $t_0\in\R$ such that for every $t\geq t_0$, $E_{\delta} \cap (E_{\delta}-\a_1t) \cap \ldots \cap (E_{\delta}-\a_kt)\neq \emptyset$.  

We use the following modified version of the correspondence principle of Furstenberg, Katznelson, and Weiss \cite{FKW}:
\begin{theorem}\label{btheorem1.5}  Suppose $E\subset \R$ with $D^*(E)>0$.  Then there exists an ergodic measure preserving flow $(X,\X,\m,\{T_t\})$ and some $\widetilde{E} \in \X$ with $\m(\widetilde{E})\geq D^*(E)$ such that if $\{ u_1,u_2,\ldots,u_k \} \subseteq \R$, then for all $\delta >0$, $$ D^*\bigr(E_{\delta} \cap (E_{\delta}-u_1) \cap \cdots \cap (E_{\delta}-u_k)\bigl) \geq \m(\widetilde{E}\cap T_{u_1}^{-1}\widetilde{E} \cap \cdots \cap T_{u_k}^{-1}\widetilde{E}).$$
\end{theorem}
The original correspondence principle\footnote{An $\R^d$ version was subsequently used by Ziegler in \cite{Z} to study configurations in $\R^d$, by examining discrete time averages for transformations which arise from an $\R^d$-action.} of Furstenberg, Katznelson, and Weiss was developed in order to study configurations in the plane and states that $E_{\delta} \cap (E_{\delta}-u_1) \cap \cdots \cap (E_{\delta}-u_k)$ is nonempty, but does not give a lower bound for the upper density, and does not guarantee that the flow $(X,\X,\m,\{T_t\})$ will be ergodic (see \cite{FKW}).  The proof of Theorem \ref{btheorem1.5} is similar to the proof in \cite{FKW}, but for the sake of completeness we include a proof in Appendix \ref{appendixCP}.  The proof also makes use of the Ergodic Decomposition Theorem and the fact that almost every point in $X$ is \textbf{\emph{quasi-generic}} (see Appendix \ref{appendixCP} for the definition) to obtain the lower bound.  

Combining Theorem \ref{btheorem1.3} and Theorem \ref{btheorem1.5} we have:
\begin{theorem}\label{btheorem1.6}  Suppose $E\subset \R$ with $D^*(E)>0$ and $\{p_1,\ldots,p_k\colon \R^d \rightarrow \R\}$ is a collection of polynomials with $p_1(0)=\ldots=p_k(0)=0$ and with complexity 0 or 1.  Then the set
\begin{equation*} \{ \bfs \in \R^d \colon\,\text{$\forall$}\,\,\delta>0, \, D^*\bigr(E_{\delta} \cap (E_{\delta}-p_1(\bfs)) \cap \cdots \cap (E_{\delta}-p_k(\bfs))\bigl)>D^*(E)^{k+1}-\varepsilon  \}  \end{equation*}
is syndetic. 
\end{theorem}  

For example, Theorem \ref{btheorem1.3} holds for the families $\{t,2t\}$, $\{t,t^2,3t^2+\pi t\}$, and $\{t,t^2+t,\ldots,t^k+t^{k-1}\}$.  It is an open question as to whether Theorem \ref{btheorem1.6} still holds when $E_{\delta}$ is replaced by $E$.    

It also follows that the conclusion of Theorem \ref{btheorem1.6} holds for a family of polynomials with complexity 2:
\begin{theorem}\label{btheorem1.7}  Suppose $E\subset \R$ with $D^*(E)>0$ and $l,m\in \N$.  Let $p \colon \R^d \rightarrow \R$ be a polynomial with $p(0)=0$ and let $\varepsilon >0$.  Then the set of $\bfs \in \R^d$ such that for all $\delta > 0$, $$D^*\bigr(E_{\delta}  \cap (E_{\delta}-m p(\bfs)) \cap ( E_{\delta} -l p(\bfs))   \cap (E_{\delta}-(l+m) p(\bfs))\bigl) >D^*(E)^4-\varepsilon $$ is syndetic. 
\end{theorem}

\subsection{Guide to the paper.}\label{bsec1.4}  We begin by giving some background information in Section \ref{bsec2}.  In Section \ref{bsec3} we show that for ergodic flows the average \eqref{maineq} is bounded by the Host-Kra seminorms, as developed in \cite{HK}, starting first with the linear case and then proving the general case using an induction argument, as developed in \cite{Berg}.  From results in \cite{HK} and \cite{Z} we then show that the Host-Kra factors are characteristic for \eqref{maineq} and hence reduce to the case where $(X,\X,\m,\{T_t\})$ is an inverse limit of nilflows.

In Section \ref{bsec4} we complete the proof of Theorem \ref{mainconvg} by reducing to the case where $(X,\X,\m,\{T_t\})$ is an ergodic nilflow.  Convergence in this setting follows from \cite{Shah}. 

In Section \ref{bsec5} we give a formula for the limit \eqref{maineq}.   First we prove Theorem \ref{btheorem1.2} using methods given in \cite{FK2}, by reducing to the case of a nilflow, then further reducing to the abelianization and using the Weyl Equidistribution Theorem.  We then show how in general the form of the limit \eqref{maineq} can be deduced from \cite{Le4}.  In particular, we show that it suffices to compute the limit of \eqref{maineq} for collections of linear polynomials.  Using this fact, we develop a method for bounding the complexity of a collection polynomials.

Section \ref{bsec6} contains the proofs of Theorems \ref{btheorem1.3} and \ref{btheorem1.4}, using techniques developed in \cite{Fran}.  The proof of Theorem \ref{btheorem1.3} makes use of the fact that the Kronecker factor is characteristic for the average \eqref{maineq} in the relevant case, allowing us to compute the limit along some syndetic set of times.  The proof of Theorem \ref{btheorem1.4} is similar, but uses the symmetry of the polynomials $\{l p, m p , (l+m) p\}$ to compensate for the fact that the characteristic factor is non-abelian.

\section{Background.}\label{bsec2} 
\subsection{The setting.}\label{bsec2.1} For simplicity of notation, we assume that all functions are real-valued, but note that all statements hold in the case of complex-valued functions.  

Throughout, $(X,\X,\m)$ is a Lebesgue space with $\m(X)=1$, and $\{T_t\}_{t\in\R}$ is a \textbf{\emph{measure preserving flow}}.  This means $\{T_t\}_{t\in\R}$ is a collection of invertible measure preserving transformations $ \{T_t \colon (X,\X,\m) \rightarrow (X,\X,\m)\}$ such that the map $\R \times X \rightarrow X$ given by $(t,x)\mapsto  T_t(x) $ is measurable, $T_0$ is the identity transformation, and $T_s \circ T_t = T_{t+s}$ for all $s,t\in\R$.  We also assume $(X,\X,\m,\brac{T_t})$ is \textbf{\emph{ergodic}}, i.e., a set $A\in\X$ satisfies $T_t(A)=A$ for all $t\in\R$ if and only if $\m(A)=0$ or $1$.  If $T \colon X \rightarrow X$ is a measure preserving transformation we frequently denote $f \circ T$ by $T f$.  

Of particular importance, as $(X,\X,\m)$ is a Lebesgue space, the map $\R \times L^2(\m) \rightarrow L^2(\m)$ given by $(t,f)\mapsto f\circ T_t $ is continuous (see \cite{AK}).  This fact allows us to work under connectedness assumptions which make several proofs simpler than the discrete counterparts, and in some cases lead to stronger results.    

We utilize the following result of Pugh and Shub. 
\begin{theorem}[Pugh and Shub, \cite{PS}]\label{btheorem2.1}  Let $\brac{T_t}_{t\in\R}$ be an ergodic measure preserving flow on a Lebesgue space $(X,\X,\m)$.  Then there exists a countable set $E\subset \R$ such that for each $t_0\notin E$, the transformation $T_{t_0}$ is ergodic.
\end{theorem}
We call $E=E(\brac{T_t})$ the \textbf{\emph{exceptional set}} of $\{T_t\}_{t\in\R}$.  


\subsection{Factors.}\label{bsec2.2}  A measure preserving flow $(Y,\Y,\nu,\{S_t\}_{t\in\R})$ is a \emph{\textbf{factor}} of the measure preserving flow $(X,\X,\m,\{T_t\}_{t\in\R})$ if there is some $\{T_t\}$-invariant, full measure subset $X'$ of $X$, some $\{S_t\}$-invariant, full measure subset $Y'$ of $Y$, and some measurable map $\pi \colon X' \rightarrow Y'$ such that $\nu=\mu \circ \pi^{-1}$ and $S_t \circ \pi (x)=\pi \circ T_t(x)$ for all $t\in \R$ and for all $x\in X'$. 

A factor $(Y,\Y,\nu,\{S_t\}_{t\in\R})$ of $(X,\X,\m,\{T_t\}_{t\in\R})$ can be naturally identified with the $\{T_t\}$-invariant sub-$\sigma$-algebra $\pi^{-1}(\Y)$ of $\X$, or equivalently, with the closed $\{T_t\}$-invariant subspace $L^2(\pi^{-1}(\Y))$ of $L^2(\X)$.  If $\Y$ is a $\{T_t\}$-invariant sub-$\sigma$-algebra of $\X$ and $f\in L^{2}(\X)$, then the \emph{\textbf{conditional expectation}} of $f$ on $\Y$ is the orthogonal projection of $f$ on the closed subspace $L^2(\pi^{-1}(\Y))$ of $L^2(\X)$, and is denoted by $\E(f|\Y)$.  

We say $(X,\X,\m,\{T_t\})$ is an \textbf{\emph{inverse limit}} of the factors 
$(X,\X_i,\m,\{T_t\})$ if $\X_i$ is an increasing sequence of $\{T_t\}_{t\in\R}$-invariant sub-$\sigma$-algebras of $\X$ and $\X= \bigvee_{i=1}^{\infty} \X_i$ up to sets of measure zero.

\subsection{Host-Kra seminorms and factors.}\label{bsec2.3}  Let $T$ be an ergodic measure preserving transformation on $(X,\X,\m)$.  In \cite{HK}, Host and Kra developed a sequence of seminorms $\{\seminorm{\cdot}_{k,T}\}_{k\in\N}$ on $L^{\infty}(\m)$ which they used to bound discrete time multiple ergodic averages.  We review some constructions and statements given in \cite{HK}. 

A collection of measure preserving systems $\{(X^{[k]},\X^{[k]},\m^{[k]},T^{[k]})\}_{k \in \N}$ is inductively defined such that $(X^{[0]},\X^{[0]},\m^{[0]})=(X,\X,\m)$, and for every integer $k \geq 1$,  $X^{[k]}=X^{2^k} $, and $T^{[k]}=T \times T \times \ldots \times T$ ($2^k$ times).  Furthermore, if $\I^{[k]}$ denotes the $T^{[k]}$-invariant $\sigma$-algebra of $(X^{[k]},\mu^{[k]},T^{[k]})$, then $\mu^{[k]}$ is defined on $X^{[k]}$ by 
$$  \int_{X^{[k]}} F\times G \, d  \mu^{[k]} = \int_{X^{[k-1]}} \E(F| \I^{[k-1]} )  \E(  G| \I^{[k-1]} )\, d  \mu^{[k-1]} $$
for all $F,G\in L^{\infty}(X^{[k-1]})$.  It follows that $\m^{[k]}$ is $T^{[k]}$-invariant.  For each $k\geq 1$ define 
$$  \seminorm{f}_k^{2^k}   \colon  =   \int_{X^{[k]}}  \bigotimes_{\varepsilon \in \{0,1\}^{{k}}} f  \, d \mu^{[k]}   $$
for all $f \in L^{\infty} (\m)$.  It was shown this defines a seminorm on $L^2(\m)$.  We sometimes write $\bigotimes_{2^{k}} f$ instead of $\bigotimes_{\varepsilon \in \{0,1\}^{{k}}} f$.

By ergodicity, the $\sigma$-algebra $\I^{[0]}$ is trivial, $\mu^{[1]}= \mu \times \mu$, and $ \seminorm{f}_1 = \abs{\int f(x) d \mu (x) } $.  Furthermore, for every integer $k \geq 1$ and every $f\in L^{\infty}(\m)$, $\seminorm{f}_{k+1}^{2^{k+1}} = \limN \avgs{N} \seminorm{f \cdot T^n f}_k^{2^k}$ and $\seminorm{f}_{k+1} \geq \seminorm{f}_k$.

Furthermore, for each ergodic measure-preserving system $(X,\X,\m,T)$, there exists a sequence of factors $\sZ_0(T) \subseteq \sZ_1 (T)\subseteq \ldots \subseteq \sZ_k(T) \subseteq \ldots$ such that for each $k\geq 1$, $\sZ_{k-1}(T)$ is \emph{characteristic} for the average $ \avgs{N}  T^{n}f_1 \cdot   T^{2n}f_2 \cdot \ldots \cdot T^{kn}f_k$ where $ f_1,\ldots,f_k \in L^\infty(\m)$.  In other words, the $L^2$-limit of this average is unchanged if $f_1,\ldots,f_k$ are replaced by $\E(f_1|\sZ_{k-1}(T)),\ldots,$ $\E(f_k|\sZ_{k-1}(T))$.  These factors are controlled by the seminorms $\brac{\seminorm{\cdot}_{k,T}}$ in the sense that for all $k\geq 1$ and for all $f\in L^{\infty}(\m)$, $\seminorm{f}_{k,T}=0$ if and only if $\E(f|\sZ_{k-1}(T))=0$.  Moreover, it is proved in \cite{HK} that each $\sZ_k(T)$ is the inverse limit of a sequence of $(k-1)$-step nilsystems.  In particular, $\sZ_0(T)$ is the trivial factor of $\X$ and $\sZ_1(T)$ is the Kronecker factor. 

\subsection{Seminorms and factors for flows.}\label{bsec2.4}  Frantzikinakis and Kra showed in \cite{FK} that if $T$ and $S$ are commuting ergodic transformations of a probability space $(X,\X,\m)$ with associated Host-Kra seminorms $\{\seminorm{\cdot}_{k,T}\}_{k\in\N}$ and $\{\seminorm{\cdot}_{k,S}\}_{k\in\N}$, then $\seminorm{f}_{k,T}=\seminorm{f}_{k,S}$ for all integers $k\geq 1$ and for all $f\in L^{\infty}(\m)$.  Furthermore, the Host-Kra factors associated to $T$ and $S$ agree.  These two facts, in combination with Theorem \ref{btheorem2.1}, allow us to define a collection of seminorms on $L^{\infty}(\m)$ corresponding to the ergodic flow $(X,\X,\m,\{T_t\})$, as well as an associated sequence of factors.
\begin{definition}\label{adef2.2} $\seminorm{f}_{k} \colon\!\!\!=\seminorm{f}_{k,T_s} $ for all $f\in L^{2}(\m)$, $s \in E^C$, and $k\in\N$.
\end{definition} 
\begin{definition}\label{adef2.3}  $\sZ_k(X,\{T_t\})\colon\!\!\!=\sZ_k(X,T_s)$ for each integer $k\geq 0$ and for all $s\in E^C$.
\end{definition}

We simply write $\sZ_k$ instead of $\sZ_k(X, \{T_t\})$ when it is clear which flow is being considered.  The use of these factors in the setting of flows was originated by Ziegler in \cite{Z}.  In the discrete setting, the $\sZ_k$ were shown to be inverse limits of nilsystems in \cite{HK}.  In \cite{Z}, Ziegler shows that the analogous result for flows holds as well.  In other words, for each integer $k\geq 0$,
 \begin{equation}\label{beq2.1}\text{$\sZ_k$ is an inverse limit of $(k-1)$-step nilflows.} \end{equation}

\section{Averages are controlled by seminorms.}\label{bsec3}
\subsection{Linear averages are controlled by seminorms.}\label{bsec3.1}  Any finite collection of polynomials $p_1,\ldots,p_k \colon \R^d \rightarrow \R$ is called a \emph{\textbf{family}}.  A family of polynomials $\{p_1,\ldots,p_k\}$ is said to be \emph{\textbf{essentially distinct}} if $p_i-p_j$ is non-constant for all $i,j\in\{1,\ldots,k\}$ with $i\neq j$ and \emph{\textbf{nice}} if the $p_i$ are non-constant and essentially distinct.  We prove:
\begin{proposition}\label{bprop3.1}  Let $\{T_t\}$ be an ergodic flow on a Lebesgue space $(X,\X,\m)$ and let $p_1,\ldots,p_k \colon \R^d \rightarrow \R$ be a nice family of linear polynomials with $p_i(0)=0$ for $i=1,\ldots,k$.  Then for all $f_1, \ldots,f_k \in L^{\infty}(\m)$ with $\norm{f_1}_{\infty},\ldots,\norm{f_k}_{\infty}\leq 1$,
$$ \limsup_{R_1,\ldots,R_d \rightarrow \infty} \norm{\avgRd T_{p_1(\bfs)} f_1 \cdot \ldots \cdot T_{p_k(\bfs)} f_k  \, d\bfs}_{L^2(\m)} \leq  \min_{1 \leq l \leq k}\seminorm{f_l}_{k}. $$
\end{proposition}

We say that a collection of transformations $\brac{T_{\a}\colon X \rightarrow X}_{\a\in\Lambda} $ is \emph{\textbf{totally ergodic}} if $T_{\a_1}^{n_1}\cdot \ldots \cdot T_{\a_l}^{n_l}$ is ergodic for all distinct elements $\a_1,\ldots,\a_l \in \Lambda$ and for all $n_1,\ldots,n_l \in \Z$ with $(n_1,\ldots,n_l) \neq (0,\ldots,0)$.  We remark that the set of zeros of a nonzero polynomial $p \colon \R^d \rightarrow \R$ has Lebesgue measure zero.  Consequently, given a nice family of polynomials $\{q_1,\ldots,q_l \colon \R^d \rightarrow \R \}$, there exists some $\Delta \in \R^d$ of Lebesgue measure zero such that $\{T_{q_1(\bfs)},\ldots,T_{q_{l}(\bfs)}\}_{\bfs\in\R \backslash \Delta}$ is a totally ergodic collection of transformations.
 
\begin{lemma}\label{blemma3.2}  For all integers $d,k\geq 1$, for each non-constant linear polynomial $p \colon \R^d \rightarrow \R$ with $p(0)=0$, and for and every $f\in L^{\infty}(\m)$,
\begin{equation}\label{beq3.1}  \seminorm{f}_{k+1}^{2^{k+1}} =  \lim_{R_1,\ldots,R_d \rightarrow \infty} \avgi{R_1} \ldots\avgi{R_d} \seminorm{f \cdot T_{p(\bfs)}f}_k^{2^k} \, d\bfs .\end{equation}
\end{lemma}

\begin{proof}  Let $\a_{1},\ldots,\a_{d}\in \R\backslash \{0\}$ such that $p(\bfs)=\a_{1}s_1+\ldots + \a_{d}s_d$ for all $\bfs =(s_1,\ldots,s_d)\in \R^d$.  For all $N_1, \ldots,N_d\in\N$ and $\bfs \in \R^d$,
\begin{eqnarray}\label{beq3.2} \lefteqn{\avgNsub{1} \ldots \avgNsub{d} \seminorm{f \cdot T_{\a_{1}}^{n_1}\ldots T_{\a_{d}}^{n_d}(T_{p(\bfs)}f)}_k^{2^k} }\\
\nonumber & = \displaystyle{ \int_{X^{[k]}} \avgNsub{1}\ldots\avgNsub{d}\Bigl[ \bigl( \bigotimes_{2^k} f \bigr) \cdot (T_{\a_{1}}^{[k]})^{n_1}\ldots (T_{\a_{d}}^{[k]})^{n_d} \bigl( \bigotimes_{2^k} T_{p(\bfs)}f\bigr)   \Bigr]  \,   d \mu^{[k]}  .}
\end{eqnarray}

First suppose $\{T_{\a_{1}},\ldots,T_{\a_{d} }\}$ is totally ergodic.  It was shown in \cite{FK} that if $T$ and $S$ are two commuting ergodic transformations of $(X,\X,\m)$, then $T^{[k]}$ and $S^{[k]}$ have the same invariant sets.  For almost every $\bfs\in\R^d$, the collection $\{T_{p(\bfs)},T_{\a_1},\ldots,T_{\a_{d}}\}$ is totally ergodic, and thus by \eqref{beq3.2}, the definition of the measures $\m^{[k]}$, the invariance of $\I^{[k]}$ under the collection $\{T_{p(\bfs)}^{[k]},T_{\a_{1}}^{[k]},\ldots,T_{\a_{d}}^{[k]} \}$, and the ergodic theorem,
\begin{eqnarray}\nonumber\lefteqn{\lim_{N_1,\ldots,N_d\rightarrow \infty} \avgNsub{1} \ldots \avgNsub{d} \seminorm{f \cdot T_{\a_{1}}^{n_1}\ldots T_{\a_{d}}^{n_d}(T_{p(\bfs)}f)}_k^{2^k}} \\
 \label{beq3.3} & & = \int  \Bigl( \bigotimes_{2^k}f \Bigr) \cdot \E(\bigotimes_{2^k}T_{p(\bfs)} f    | \I^{[k]}) \,d \mu^{[k]}  =    \seminorm{f}_{k+1}^{2^{k+1}}  .
\end{eqnarray}

For $R\in \R$, let $\floor{R}$ denote the integer part of $R$.  As the integrand of \eqref{beq3.1} is bounded, for each $i\in\{1,\ldots,d\}$, we can replace $R_i$ with $\floor{R_i}$ without changing the limit.  Furthermore, we can write $[0,\floor{R_i}-1]= \bigcup_{n_i=1}^{\floor{R_i}-1} [n_i,n_i+1]$ for each $i\in\{1,\ldots,d\}$ and break up the integrals accordingly.  Thus by \eqref{beq3.3} and the linearity of $p$,
 \begin{eqnarray*}\lefteqn{ \lim_{R_1,\ldots,R_d \rightarrow \infty} \frac{1}{\prod_{i=1}^{d}R_i} \int_{0}^{R_1}\ldots \int_{0}^{R_d}   \seminorm{f\cdot T_{p( \bfs)}f}_k^{2^k} \, d\bfs }\\
  & = &\lim_{R_1,\ldots,R_d \rightarrow \infty} \frac{1}{\prod_{i=1}^{d} \floor{R_i}}\sum_{n_1=0}^{\floor{R_1}-1} \int_{n_1}^{(n+1)} \ldots \sum_{n_d=0}^{\floor{R_d}-1} \int_{n_d}^{(n+1)}   \seminorm{f\cdot T_{p( \bfs)}f}_k^{2^k} \, d\bfs \\
 & = & \int_{[0,1]^d}  \lim_{R_1,\ldots,R_d \rightarrow \infty}\frac{1}{\prod_{i=1}^{d}\floor{R_i}}\sum_{n_1=0}^{\floor{R_1}-1} \ldots \sum_{n_d=0}^{\floor{R_d}-1}\seminorm{f \cdot T_{\a_{1}}^{n_1}\ldots T_{\a_{d}}^{n_d}(T_{p ( \bfs)} f)}_k^{2^k} \, d\bfs \\
 & = &   \int_{[0,1]^d} \seminorm{f}_{k+1}^{2^{k+1}} d\bfs = \seminorm{f}_{k+1}^{2^{k+1}}. 
\end{eqnarray*}

If $\{T_{a_{1}},\ldots,T_{a_{d} }\}$ is not totally ergodic, fix $u \in \R$ such that $u>0$ and $\{T_{a_{1}u},\ldots,T_{a_{d}u }\}$ is totally ergodic.  Let $\{\widetilde{T}_t\}_{t\in\R}$ be the flow given by $\widetilde{T}_t= T_{ut}$ for all $t\in\R$.  Then $\{\widetilde{T}_{a_{1}},\ldots,\widetilde{T}_{a_{d} }\}$ is totally ergodic and hence \eqref{beq3.1} holds when $T_t$ is replaced with $\widetilde{T}_t$.  The change of variable $(s_1,\ldots,s_d) \mapsto (us_1,\ldots,us_d)$ now gives the result.  
\end{proof}

We now prove Proposition \ref{bprop3.1} using a version of the van der Corput Lemma and a corollary.  For a full statement and proof, see Lemma \ref{vdc} and Corollary \ref{bcorB.2}, Appendix \ref{appendixVDC}.  The use of van der Corput's Lemma for bounding discrete time averages was first introduced by Bergelson in \cite{Berg}. 
\begin{proof}[Proof of Proposition \ref{bprop3.1}]  We proceed by induction on $k$.  First suppose $k=1$.  For $\bfs\in\R^d$ we apply the van der Corput Lemma to the elements $g_{\bfs}=T_{p_1(\bfs)}f_1$ in $L^2(\m)$.  For any set $\Psi \subseteq \R^d$ of finite positive Lebesgue measure,
\begin{align}\label{beq3.4} & \displaystyle{\limsup_{N\rightarrow \infty} \norm{  \avgRd T_{p_1(\bfs)} f_1 \, d\bfs }_{L^2(\m)}^2} \\
\nonumber & \leq  \displaystyle{\limsupN \frac{1}{m(\Psi)^2} \int_{\Psi} \int_{\Psi} \avgRd \int T_{p_1(\bfs+\bfu)} f_1\cdot  T_{p_1(\bfs+\bfv)} f_1 \, d\m \, d\bfs \, d\bfu \, d\bfv}\\
\nonumber & =  \displaystyle{\limsupN \frac{1}{m(\Psi)^2} \int_{\Psi} \int_{\Psi} \avgRd \int T_{p_1(\bfu)} f_1\cdot  T_{p_1(\bfv)} f_1 \, d\m \, d\bfs \, d\bfu \, d\bfv}\\
\nonumber  & = \displaystyle{ \int \frac{1}{m(\Psi)^2} \int_{\Psi} \int_{\Psi}T_{p_1(\bfu)} f_1 \cdot T_{p_1(\bfv)}f_1  \, d\bfu \, d\bfv \, d\m  }.\end{align}
By taking the $\limsup$ over all rectangles ${\Psi} \subset \R^d$ and by the ergodic theorem, we see that \eqref{beq3.4} is less than or equal to $ \abs{\int f_1 \,d\m}^2 =  \seminorm{f_1}_1^2$.

Next suppose $k\geq 2$ and Proposition \ref{bprop3.1} holds for $k-1$.  We show Proposition \ref{bprop3.1} also holds for $k$.  For $\bfs\in \R^d$ we apply the van der Corput Lemma and Corollary \ref{bcorB.2} to the element $g_\bfs=T_{p_1(\bfs)} f_1 \cdot \ldots \cdot T_{p_k(\bfs)} f_k$ of $L^2(\m)$. For any ${\Psi}\subset \R^d$ with positive finite Lebesgue measure and for any $l\in\{1,\ldots,k-1\}$ (the case $k=l$ is similar), 
\begin{eqnarray}\label{beq3.5}\lefteqn{\displaystyle{ \limsup_{N \rightarrow \infty} \norm{ \avgRd  \prod_{i=1}^{k} T_{p_i(\bfs)} f_i \ d\bfs }_{L^2(\m)}^2}}\\ 
\nonumber & \leq & \displaystyle{\frac{1}{m({\Psi})} \int_{\Psi} \frac{1}{m({\Psi})} \int_{\Psi} \limsup_{N \rightarrow \infty}\norm{T_{p_k(\bfu)}f_k\cdot T_{p_k(\bfv)}f_k}_{L^2(\m)} \cdot}\\
 \nonumber & & \displaystyle{\norm{\avgRd    \prod_{i=1}^{k-1} T_{p_i(\bfs)-p_k ( \bfs) } \bigl( f_i\circ T_{p_i(\bfu)}\cdot  f_i\circ T_{p_i ( \bfv)}\bigr)   \, d\bfs }_{L^2(\m)} \ d\bfu\,d\bfv } \\
 \nonumber &\leq & \displaystyle{\bigl(\frac{1}{m({\Psi})} \int_{\Psi}  \frac{1}{m({\Psi})} \int_{\Psi}\seminorm{f_l\cdot T_{p_l ( \bfv)-p_1(\bfu)} f_l}_{k-1}^{2^{k-1}} \, d\bfu \,d\bfv \bigr)^{\frac{1}{2^{k-1}}}} .
 \end{eqnarray}
 
Notice that the map $(\bfu,\bfv) \mapsto p_l(\bfv)-p_1(\bfu)$ is a linear polynomial from $\R^{2d}$ into $\R$.  By taking the $\limsup$ over all rectangles ${\Psi} \subset \R^d$ and using Lemma \ref{blemma3.2}, we see that \eqref{beq3.5} is less than or equal to $ \seminorm{f_l}_{k}^2$.
\end{proof}

\subsection{Polynomial averages are controlled by seminorms.}\label{bsec3.2}  In this section we prove the following:
\begin{proposition} \label{bprop3.3}  Let $\{T_t\}$ be an ergodic flow on a Lebesgue space $(X,\X,\m)$.  For any $k \in \N$ and for any nice family of polynomials $P=\{p_1,\ldots, p_k \colon \R^d \rightarrow \R\}$ with $p_i(0)=0$ for $i=1,\ldots,k$, there exists $r \in \N$ such that for any $f_1, \ldots, f_k \in L^\infty(\m)$,
$$ \limsup_{R_1,\ldots,R_d \rightarrow \infty} \norm{ \avgRd T_{p_1(\bfs)}f_1 \cdot \ldots \cdot T_{p_k(\bfs)}f_k \ d\bfs}_{L^2(\m)}\leq \min_{1 \leq l \leq k}\seminorm{f_l}_{r}.$$
\end{proposition}
\begin{remark}  The integer $r$ in Proposition \ref{bprop3.3} depends neither on the flow $(X,\X,\m,\brac{T_t})$ nor on $d$.
\end{remark}

The following is a consequence of Propositions \ref{bprop3.1} and \ref{bprop3.3}:
\begin{corollary} \label{bcor3.4}  Let $\brac{T_t}$ be an ergodic flow on a Lebesgue space $(X,\X,\m)$.  For any nice family of polynomials $\{p_1,\ldots, p_k \colon \R^d \rightarrow \R\}$, there exists $r\in\N$ such that for all $f_1, \ldots, f_k \in L^\infty(\m)$,  
$$ \norm{ \avgRd \prod_{i=1}^{k} T_{p_i(\bfs)}f_i  \,d\bfs -\avgRd \prod_{i=1}^{k}T_{p_i(\bfs)}\E(f_i |\sZ_{r}) \,d\bfs         }_{L^2(\m)} $$
converges to zero as $R_1,\ldots,R_d \rightarrow \infty$.  If $\{p_1,\ldots, p_k \}$ are all linear then $r=k-1$.
\end{corollary}

In other words, Corollary \ref{bcor3.4} states that $\sZ_{r}$ is \emph{\textbf{characteristic}} for the average
\eqref{maineq}.  Leibman proved the discrete time version of Corollary \ref{bcor3.4} in \cite{L1}; our proof (including elements of the proof of Proposition \ref{bprop3.3}) is similar.  
 \begin{proof}[Proof of Corollary \ref{bcor3.4}]  By the multilinearity of the average it suffices to show that
\begin{equation}\label{beq3.6} \lim_{R_1,\ldots,R_d \rightarrow \infty}  \norm{ \avgRd \prod_{i=1}^{k} T_{p_i(\bfs)}f_i  \,d\bfs }_{L^2(\m)}  =0 \end{equation}
 whenever $\E(f_i|\sZ_{r})=0$ for some $i\in\{1,2,\ldots,k  \}$.  Notice that $\E(f_i|\sZ_{r})=0$ exactly when $\E(T_{p_i(0)}f_i|\sZ_{r})= T_{p_i(0)}\E(f_i|\sZ_{r})=0$.  It follows from definitions \ref{adef2.2} and \ref{adef2.3} that $\E(f_i|\sZ_{r})=0$ if and only if $\seminorm{T_{p_i(0)}f_i}_{r+1}=0$ and hence \eqref{beq3.6} follows from Propositions \ref{bprop3.1} and \ref{bprop3.3}.
 \end{proof} 

We prove Proposition \ref{bprop3.3} using an induction argument, as developed by Bergelson in \cite{Berg}.  If $P=\brac{p_1,\ldots,p_k}$ is a family of polynomials then its \emph{\textbf{degree}}, $\deg P$, is the largest degree of its elements.  We define two polynomials $p$ and $q$ to be \emph{\textbf{equivalent}} if $\deg p=\deg q$ and $\deg\abs{p-q}< \deg p$.  For example, $t^2+t$ and $t^2$ are equivalent, while $t^2+t$ and $3t^2$ are not.  This partitions the set of all polynomials into equivalence classes, and the degree of an equivalence class is the degree of any of its elements.  

We assign each family $P$ of degree $b$ a \emph{\textbf{weight vector}} $\om(P)=(\om_1,\ldots,$\\$\om_b)\in \N^b$, where each $\om_i$ is the number of equivalence classes of degree $i$ in $P$, and we say $\om(P)$ has degree $b$.  For example, the weight vector of $\{t,2t,3t, t^2,t^2-t,4t^2+t,t^3  \}$ is $(3,2,1)$.  We write $\om < \om'$ if $\deg \om < \deg \om' $.  If $\deg \om = \deg \om'$, we resort to right-aligned lexicographical ordering.  In other words, $\om < \om'$ if $\deg \om < \deg \om' $, or if $\deg \om = \deg \om'$ and there exists some $j\leq b$ so that $\om_{j}<\om_{j}'$ and $\om_{i}=\om_{i}'$ for $j<i\leq b$.  The set of weight vectors is well ordered with respect to this relation, and we use induction on this set.

We call a nice family of polynomials $P=\{p_1,\ldots,p_k\}$ \emph{\textbf{standard}} if $\deg P= \deg p_1$. 

\begin{proof}[Proof of Proposition \ref{bprop3.3}]  We first prove that for every standard family $P=\{p_1,\ldots, p_k \colon \R^d \rightarrow \R\}$ with $p_i(0)=0$ for $i=1,\ldots,k$, there exists $r \in \N$ such that for any $f_1, \ldots, f_k \in L^\infty(\m)$,
\begin{equation}\label{beq3.7} \limsup_{R_1,\ldots,R_d \rightarrow \infty} \norm{ \avgRd T_{p_1(\bfs)}f_1 \cdot \ldots \cdot T_{p_k(\bfs)}f_k \ d\bfs}_{L^2(\m)}\leq \seminorm{f_1}_{r}.\end{equation}

We proceed by induction on $\om=\om(P)$.  Proposition \ref{bprop3.1} is the base case in our induction.  Let $P= \brac{p_1, \ldots,p_k \colon \R^d \rightarrow \R}$ be a standard family of degree $\geq 2$ and of weight $\om$, and suppose that \eqref{beq3.7} holds for any standard family with weight vector $\om'<\om$.  We assume that $p_k$ is a polynomial of minimal degree in $P$.  Without loss of generality, we assume that $\norm{f_1}_{\infty},\ldots,\norm{f_k}_{\infty} \leq 1$.  Let $I_1 = \brac{i \in \brac{1, \ldots,k} \colon \deg p_i=1}$ and $I_2 = \brac{i \in \brac{1, \ldots,k} \colon \deg p_i\geq 2}$.

We use the van der Corput Lemma and Corollary \ref{bcorB.2}.  Write $ g_{\bfs}(x)=T_{p_1(\bfs)}f_1 \cdot \ldots \cdot T_{p_k(\bfs)}f_k $ for every $\bfs\in\R^d$.  Then
\begin{eqnarray*}\nonumber\lefteqn{\displaystyle{ \avgRd \left\langle  g_{\bfs+\bfu},g_{\bfs+\bfv} \right\rangle \,d\bfs } }\\
\nonumber & = & \displaystyle{\avgRd\int\prod_{i \in I_2} T_{p_i(\bfs+\bfu)}f_i\cdot \prod_{i \in I_2} T_{p_i(\bfs+\bfv)}f_i}\\
\nonumber & & \displaystyle{\hspace{2in} \cdot  \prod_{i \in I_1} T_{p_i(\bfs+\bfv)} \bigl(f_i \cdot T_{p_i(\bfu)-p_i(\bfv) } f_i  \bigr)
 \,  d\m\,d\bfs}\\
 & =& \displaystyle{ \avgRd\int \prod_{j=1}^m T_{q_{\bfu,\bfv,j}(\bfs)} h_{\bfu,\bfv,j} 
 \, d\m \,d\bfs} \end{eqnarray*}
where, for $\bfu,\bfv \in \R^d$, $q_{\bfu,\bfv,1}, \ldots, q_{\bfu,\bfv,m}$ are the elements of the family 
 $$P_{\bfu,\bfv}= \brac{p_i(\bfs+\bfu), \ p_i(\bfs +\bfv) \colon i \in I_2} \bigcup \brac{p_i(\bfs+\bfv) \colon i \in I_1} , $$ and each $h_{\bfu,\bfv,j}$ is of the form $f_i$ for some $i \in I_2$, or $f_i \cdot  T_{p_i(\bfu)-p_i(\bfv) }f_i  $ for some $i \in I_1$.  We assume that $q_{\bfu,\bfv,1}(\bfs)=p_{1}(\bfs+\bfv)$ and $q_{\bfu,\bfv,m}(\bfs)=p_{k}(\bfs+\bfv)$.

As $\{T_t\}_{t\in\R}$ is $\m$-preserving, by the Cauchy-Schwarz Inequality,
\begin{align}\nonumber\lefteqn{\avgRd \left\langle  g_{\bfs+\bfu},g_{\bfs+\bfv} \right\rangle \,d\bfs }\\
\label{beq3.8} & \leq  &\norm{h_{\bfu,\bfv,m}}_{L^2(\m)} \cdot \norm{  \avgRd    \prod_{j=1}^{m-1}T_{(q_{\bfu,\bfv,j}-q_{\bfu,\bfv,m})(\bfs)} h_{\bfu,\bfv,j} \,   d\bfs }_{L^2(\m)} .
\end{align}
For almost all $(\bfu,\bfv)\in \R^{2d}$, the collection of polynomials
$$ P_{\bfu,\bfv}'=\brac{q_{\bfu,\bfv,1}-q_{\bfu,\bfv,m}, \ldots, q_{\bfu,\bfv,m-1} -q_{\bfu,\bfv,m}}  $$
is a standard family.  Furthermore, $P$, $P_{\bfu,\bfv}$ and $P_{\bfu,\bfv}'$ have the same equivalence classes, of the same degrees, with the exception that in $P_{\bfu,\bfv}'$ the equivalence class in $P_{\bfu,\bfv}$ containing $q_{\bfu,\bfv,m}$ either splits into one or more equivalence classes of lower degree or vanishes completely.  Thus, for all $(\bfu,\bfv)\in \R^{2d}$, $\om(P_{\bfu,\bfv}')<\om(P)=\om$.  

There are only finitely many integer vectors with $\om' <\om$ which are the weights of families with $m < 2k$ elements.  Thus there exists $r \in \N$ such that for all standard families $\{Q_1, \ldots , Q_m \colon \R^d \rightarrow \R\}$ of weight $\om' < \om$ with $m \leq 2k$, and any $H_1,\ldots, H_m \in L^\infty(\m)$,
\begin{equation}\label{beq3.9} \limsup_{R_1,\ldots,R_d \rightarrow \infty}\norm{ \avgRd T_{Q_1(\bfs)}H_1 \cdot \ldots \cdot T_{Q_m(\bfs)}H_m \ d\bfs}_{L^2(\m)} \leq \seminorm{H_1}_{r}. \end{equation} 
Combining \eqref{beq3.9} and \eqref{beq3.8} and using Corollary \ref{bcorB.2},
\begin{eqnarray*}\nonumber\lefteqn{ \limsup_{R_1,\ldots,R_d \rightarrow \infty} \norm{ \avgRd T_{p_1(\bfs)}f_1 \cdot \ldots \cdot T_{p_k(\bfs)}f_k \, d\bfs}_{L^2(\m)}^2 }  \\
& \leq &  \frac{1}{m(\Psi)^2}\int_{\Psi} \int_{\Psi}  \seminorm{h_{\bfu,\bfv,1}}_r   \, d\bfu \, d\bfv = \seminorm{f_1}_{r}. 
\end{eqnarray*}

We now prove the theorem in general, where $P=\brac{p_1, \ldots, p_k}$ is a nice, but not necessarily standard, family of polynomials of degree $b$.  Let $f_1,\ldots,f_k \in L^\infty(\m)$.  By Corollary \ref{bcorB.3} there exists a F\o lner sequence $\brac{\Theta_N}_{N\in \N}$ in $\R^{3d}$ such that 
\begin{eqnarray}\label{beq3.10}\lefteqn{ \limsup_{N \rightarrow \infty}   \norm{  \avgRd \prod_{i=1}^k   T_{p_i(\bfs)}f_i \,    d\bfs        }_{L^2(\m)}^{2} }\\
\nonumber & \leq & \limsup_{N \rightarrow \infty}    \frac{1}{m(\Theta_N)} \int_{\bfu,\bfv,\bfs \in \Theta_N} \int \prod_{i=1}^k   T_{p_i(\bfs+\bfu)}f_i \cdot  \prod_{i=1}^k   T_{p_i(\bfs+\bfv)}f_i \ d\m \,   d\bfs\,d\bfu\,d\bfv\\
\nonumber & \leq & \limsup_{N \rightarrow \infty}   \Big\|  \frac{1}{m(\Theta_N)} \int_{\bfu,\bfv,\bfs \in \Theta_N}  \prod_{i=1}^k   T_{p_i(\bfs+\bfu)+q(\bfs)}f_i \cdot \\
\nonumber & &\hspace{2.5in}  \prod_{i=1}^k   T_{p_i(\bfs+\bfv)+q(\bfs)}f_i \,    d\bfs \,d\bfu\,d\bfv     \Big\|_{L^2(\m)} 
\end{eqnarray}    
for any polynomial $q\colon \R^d \rightarrow \R$ of degree $b$.  The set 
$$  \brac{p_i(\bfs+\bfu)+q(\bfs),p_i(\bfs+\bfv)+q(\bfs)\colon 1 \leq i \leq k }  $$ 
of polynomials $\R^{3d} \rightarrow \R$ is a standard family of degree $b$ with $2k$ elements.  Thus there exists $r \in \N$ such that \eqref{beq3.10} is less than or equal to $\seminorm{f_l}_r$ for each $l=1,\ldots,k$.
\end{proof}

\section{Convergence on a nilsystem.}\label{bsec4}
\subsection{Nilflows.}\label{bsec4.1}  Let G be a group.  For $h, g \in G$ we write $[g,h]=g^{-1}h^{-1}gh$.  For $A, B \subseteq G$, $[A,B]$ is the closed subgroup of $G$ spanned by $\{[a,b]$: $a \in A$, $b \in B\}$.  The \textbf{\emph{lower central series}} $G=G_1 \supset G_2 \supset \cdots \supset G_j \supset G_{j+1} \supset \cdots$ of G is defined by $G_1=G$ and $G_{j+1}=[G,G_j]$ for $j \geq 1$.  We say G is \emph{\textbf{r-step nilpotent}} if $r$ is the smallest integer such that $G_{r+1}=\{ \textit{Id}\}$.

Let $G$ be an $r$-step nilpotent Lie group and let $\Gamma$ be a \textbf{\emph{uniform}} subgroup (i.e. $\G$ is a discrete cocompact subgroup).  The compact manifold $X=G/\G$ is called an \emph{\textbf{r-step nilmanifold}}.  Let $a$ be a fixed element of $G$ and let $T_a:X \rightarrow X$ be the transformation defined by $T_a(g\G)= (a \cdot g)\G$ for all $g\in G$.  Let $\m$ be Haar measure on $X$.  Then $(X,\mu,T_a)$ is called an \emph{\textbf{r-step nilsystem}} and $T_a$ is called a \textbf{\emph{nilrotation}}.  If $\{a_t\}_{t\in\R}$ is a one-parameter subgroup of $G$ then $\{a_t\}_{t\in\R}$ induces a flow $\{T_{a_t}\}_{t\in\R}$ on $X$ defined by $T_{a_t}(g\G)= (a_t \cdot g)\G$ for all $g\in G$ and for all $t\in\R$.  A flow defined in this manner is called a \emph{\textbf{nilflow}}.

A \textbf{\emph{sub-nilmanifold}} of $X$ is a closed subset $Y$ of $X$ of the form $Y=Hx$, where $x$ is an element of $ X$ and $H$ is a closed subgroup of $G$.  If $H$ is a closed subgroup of $G$, then $H\G/\G$ is a subnilmanifold of $X$ if and only if $H \cap \G$ is uniform in $H$ if and only if $H\G$ is closed in $G$ (see \cite{L2}).  Our goal is to describe the orbits of certain paths in $X$. 

\subsection{Polynomial paths.}\label{bsec4.2}  As we only consider continuous ergodic flows, it suffices to assume $X$ is connected.  A \textbf{\emph{(multi-parameter) path}} $\{g_{\bfs}\}_{\bfs\in\R^d}$ in $G$ is a continuous function $g\colon \R^d \rightarrow G$ and we write $g_{\bfs}=g(\bfs)$ for $\bfs\in\R^d$.  If $g \colon \R^d \rightarrow G$ is a continuous homomorphism, then $g(\bfs)$ is called a \textbf{\emph{linear path}}.  Any path in $G$ naturally induces a path in $X$.  

Let $G^0$ be the connected component of the identity element in $G$.  Then $G^0\G$ is both open and closed in $G$, hence $G^0\G /\G $ is both open and closed in $X$, and $X = G^0\G /\G $.  Let $\Theta \colon G^0\G /\G \rightarrow G^0/(\G \cap G^0)$ be the map given by $\Theta(g_0 \g \G )= g_0 \G \cap G^0$ for $g_0 \in G^0$ and $\g \in \G$.  This map is a homeomorphism which preserves the left action of $G^0$.  If $g(0) \in G^0$ then $g(\bfs) \subseteq G^0$ for all $\bfs \in \R^d$, and $\Theta$ preserves the orbits of $\{g(\bfs)\}_{\bfs\in\R^d}$.  Thus, in order to describe the orbits of $\{g(\bfs)\}_{\bfs\in\R^d}$ in $X$, it suffices to describe the orbits of $\{g(\bfs)\}_{\bfs\in\R^d}$ in $G^0/(\G \cap G^0)$.  We frequently use the map $\Theta$ to reduce to the case when $G$ is connected.

If $G$ is connected then the exponential map from the Lie algebra of $G$ into $G$ is onto.  In particular, for every element $a$ in $G$ there exists some one-parameter subgroup $\{\a(t)\}_{t\in\R}$ such that $\a(1)=a$.  We denote $\a(t)$ by $a^t$.  

By \cite{M}, if $G$ is any connected simply-connected nilpotent Lie group, and $\G$ is a closed uniform subgroup of $G$, then $G$ contains a \textbf{\emph{Malcev basis}}.  In other words, there is a finite collection $\{a_1,\ldots,a_l\}\subseteq \G$ so that each $a\in G$ is uniquely representable in the form $a=a_1^{t_1}\ldots a_l^{t_l}$ for some $t_1,\ldots,t_l\in\R$.  Furthermore, every one-parameter subgroup $\{a_t\}_{t\in\R}$ of $G$ is polynomial in $\{a_1,\ldots,a_l\}$.  This means there exist polynomials $q_1 ,\ldots,q_l \colon \R \rightarrow \R$ so that $a_t=a_1^{q_1(t)} \ldots a_l^{q_l(t)}$ for all $t\in \R$.  Every connected nilpotent Lie group is a factor of a connected simply-connected nilpotent Lie group, and hence also has these properties.  Thus we may restrict our attention to \textbf{\emph{(multi-parameter) polynomial paths}}, i.e., multi-parameter paths of the form $g(\bfs)=a_1^{p_1(\bfs)}\cdot \ldots \cdot a_l^{p_l(\bfs)}$ for some $a_1,\ldots,a_l\in G$, some collection of polynomials $\{p_1,\ldots,p_l \colon \R^d \rightarrow \R\}$, and for all $\bfs \in \R^d$.   

\subsection{Uniform distribution on a subnilmanifold.}\label{bsec4.3}  A multi-parameter path $\brac{x_\bfs}_{\bfs\in\R^d}$ in $X$ is \textbf{\emph{uniformly distributed}} in $X$ if 
$$\limN \frac{m\bigl(\brac{\bfs \in \R^d \colon  x_\bfs\in U}\cap [0,R_1]\times \ldots \times [0,R_d]\bigr)}{R_1\cdot \ldots \cdot R_d}=\m(U)$$
for any open set $U$ in $X$.  Equivalently, for any $f\in C(X)$,
$$\lim_{R_1,\ldots,R_d \rightarrow \infty} \avgRd f(x_\bfs) \,d\bfs=\int f \,d\m  .$$

We use the following specific case of a more general result of Shah:

\begin{proposition}[Shah, \cite{Shah}]\label{bprop4.1}  Suppose $G$ is a nilpotent Lie group and $\G \subset G$ is a uniform subgroup.  Let $g\colon \R^d \rightarrow G$ be a polynomial path and let $x\in X=G/\G$.  Then there exists a connected closed subgroup $H$ of $G$ such that $Y=Hx$ is a closed sub-nilmanifold of $X$, $\overline{\brac{g(\bfs)x}}_{\bfs\in\R^d}=Hx$, and $\brac{g(\bfs)x}_{\bfs\in\R^d}$ is uniformly distributed in $Hx$.
\end{proposition}

In \cite{Shah}, Shah proves a more general version of Proposition \ref{bprop4.1} for real algebraic groups.  Every nilpotent Lie group is isomorphic to a real algebraic group \cite{Helg}, and hence Proposition \ref{bprop4.1} follows.  Proposition \ref{bprop4.1} follows from \cite{R} when $\{g(t)\}_{t\in\R}$ is linear.  An ergodic proof of the case where $d=1$, $\Phi_N=[0,N]$ for all $N \in \N$, $G$ is connected, and $g$ is linear is given by Green in \cite{AGH}.  Leibman proved analogous versions of Proposition \ref{bprop4.1}, as well as Corollary \ref{bcor4.2} and Proposition \ref{bprop4.4} below, for polynomial mappings from $\Z^d$ to $G$ in \cite{L3}.

\begin{corollary}\label{bcor4.2} Suppose $g\colon \R^d \rightarrow G$ is a polynomial path.  Let $x$ be any element of $X$, and let $Y=\overline{\brac{g(\bfs)x}}_{\bfs\in\R^d}$.  For any $f\in C(X)$,
$$\lim_{R_1,\ldots,R_d} \avgRd f(g(\bfs)x) \ d\bfs  =\int_Y f \ d\m_Y,$$
where $\m_Y$ is Haar measure on $Y$.
\end{corollary}

\subsection{Proof of Theorem \ref{mainconvg}.}\label{bsec4.7}
We now have all the tools necessary to prove Theorem \ref{mainconvg}.
\begin{proof}[Proof of Theorem \ref{mainconvg}]  We may always write the average so that $\{p_1,\ldots,p_k\}$ are essentially distinct.  By using the ergodic decomposition of the measure $\m$, it suffices to assume $\{T_t\}$ is ergodic.  By Corollary \ref{bcor3.4}, $\sZ_{r}$ is characteristic for the average \eqref{maineq} for some $r\in\N$, and hence it suffices to assume $\X$ is equal to $\sZ_{r}$.  By \eqref{beq2.1}, $\sZ_{r}$ is an inverse limit of $(r-1)$-step nilflows, and by an approximation argument it further suffices to assume $(X,\X,\m,\{T_t\})$ is a $(r-1)$-step nilflow.  Suppose $T_t = T_{a_t}$ for some one-parameter subgroup $\{a_t \} \subseteq G$.  We now obtain Theorem \ref{mainconvg} from Corollary \ref{bcor4.2} as follows.  Replace $X=G/\G$ with $X^k=G^k/\G^k$, $g(\bfs)$ with $(a_{p_1(\bfs)},\ldots , a_{p_k(\bfs)})$, and $f$ with $f_1 \otimes \ldots \otimes f_k$.  Applying Corollary 4.2 to points on the diagonal of $X^k$, we obtain pointwise convergence of the average \eqref{maineq} when $(X,\X,\m,\{T_t\})$ is a nilflow.  Convergence in $L^2(\m)$ for the general case follows.
\end{proof}

\subsection{Tools for computing the limit.}\label{bsec4.4}  We now give an important result that is useful for computing the limit of \eqref{maineq} in the next section.  We denote the connected component of the identity of $G$ as $G^0$.  Let $Z$ be the \textbf{\emph{maximal factor torus}} of $X$, $Z=G / ([G^0,G^0]\G)$, and let $\rho \colon X \rightarrow Z$ be the factorization mapping.  We show that well distribution on $X$ is equivalent to well distribution on $Z$.

\begin{proposition}\label{bprop4.4}  Suppose $X$ is connected, $x\in X$, and $g\colon \R^d \rightarrow G$ is a polynomial path.  The following are equivalent:
\begin{enumerate}
\item  $\brac{g(\bfs)x}_{\bfs\in\R^d}$ is dense in $X$;
\item  $\brac{g(\bfs)x}_{\bfs\in\R^d}$ is uniformly distributed in $X$;
\item  $\brac{g(\bfs)\rho(x)}_{\bfs\in\R^d}$ is dense/uniformly distributed in $Z$.
\end{enumerate}
\end{proposition}

In the case where $G$ is connected and $g$ is given by a one-parameter subgroup of $G$, Proposition \ref{bprop4.4} was shown by Green (see also \cite{P}): 
 
\begin{theorem}[Green, \cite{AGH}]\label{btheorem4.5}  If $(X=G / \G,\X,\m,\{T_t\})$ is nilflow with $G$ connected, then $\{T_t\}$ is ergodic on $X$ if and only if it is ergodic on $G/G_2\G$.
\end{theorem}
 
\begin{proof}[Proof of Proposition \ref{bprop4.4}:]  The proof is similar to the proof of Theorem B in \cite{L3}, but we state it here for the sake of completeness.

 $(1)$ implies $(2)$ by Proposition \ref{bprop4.1}.  That $(2)$ implies $(1)$ follows from the definition of \emph{well distribution} and the fact that in a compact metric space every open set has positive measure.  It is clear that (1) implies (3).

Assume (3) holds.  \emph{Case 1:}  Suppose $G$ is connected.  Then $Z=G/G_2 \G$.  By Proposition \ref{bprop4.1}, there is a closed subgroup $H$ of $G$ so that $\overline{\{g(\bfs)x\}}_{\bfs\in\R^d}=Hx$.  Therefore $Z=H \rho(x)$ and hence $G=HG_2\G$.  As $\G$ is countable, by the Baire Category Theorem $HG_2$ has non-empty interior.  Since $G$ is connected, we have $G=HG_2$.  By Lemma 3.4 in \cite{L2}, $H=G$ and thus $\overline{\brac{g(\bfs)x}}_{\bfs\in\R^d}=X$.

\emph{Case 2:}  Now assume $G$ is not necessarily connected.  Without loss of generality, we may assume $g(0)=1_G$.  Then $g(\bfs)\in G^0$ for all $\bfs\in \R^d$.  Let $\Theta \colon X \rightarrow G^0/(\G \cap G^0)$ be as defined in Section \ref{bsec4.2}.  As $\Theta$ preserves the action of $g(\bfs)$, and as $\Theta([G^0,G^0]\G)=[G^0,G^0](\G \cap G^0)$, we have that $\{g(\bfs)\Theta(x)\}_{\bfs\in\R^d}$ is well distributed in $G^0/[G^0,G^0](\G \cap G^0)$.  By case 1, $g(\bfs)$ is well distributed in $G^0/(\G \cap G^0)$, and since $\Theta$ is a homeomorphism, $g(\bfs)$ is well distributed in $ G^0\G /\G =X$.
\end{proof}

\section{Computation of the limit.}\label{bsec5}
\subsection{Independent polynomial averages converge to the product of the integrals.}\label{bsec5.1}  In this subsection we prove Theorem \ref{btheorem1.2}.  The idea of the proof is similar to, but also simpler than, that of the discrete time version given in \cite{FK2}.

We call a family of polynomials $\brac{p_1,\ldots,p_k\colon \R^d \rightarrow \R}$ \textbf{\emph{$\R$-independent}} if there does not exist a set of real numbers $\{a_1,\ldots , a_k \}$, which are not all zero, such that $a_1p_1+ \ldots + a_kp_k$ is a constant polynomial.  

By Corollary \ref{bcor3.4}, \eqref{beq2.1}, and an approximation argument, it suffices to prove the following: 

\begin{proposition}\label{bprop5.1}  Let $(X=G/\G, \sG/\G, \m , \brac{T_t})$ be an ergodic nilflow induced by a one-parameter subgroup $\{a_t\}_{t\in\R}$ of $G$.  If $\brac{p_1,p_2,\ldots,p_k\colon \R^d \rightarrow \R}$ is an $\R$-independent family of polynomials, then for every $x\in X$, the path $\{(a_{p_1(\bfs)}x,a_{p_2(\bfs)}x,\ldots,a_{p_k(\bfs)}x)\}_{\bfs\in\R^d}$ is uniformly distributed in $X^k$.
\end{proposition}

\begin{proof}  By Proposition \ref{bprop4.4} it suffices to prove Proposition \ref{bprop5.1} under the assumption that $G$ is abelian.  Let $\Theta \colon X  \rightarrow G^0/\G \cap G^0 $ be as defined in Section \ref{bsec4.2}, and let $\Theta_k=\Theta \times \ldots \times \Theta$ ($k$-times).  As the homeomorphism $\Theta_k \colon X^k  \rightarrow (G^0)^k/(\G \cap G^0)^k $ preserves the action of $\{(a_{p_1(\bfs)},a_{p_2(\bfs)},\ldots,a_{p_k(\bfs)})\}_{\bfs\in\R^d}$, we may further that assume $G$ is connected. 

As $G$ is abelian, $\G$ is a normal subgroup of $G$.  Thus $G/\G$ is a connected compact abelian Lie group and is isomorphic to some finite dimensional torus $\Tor^m$.  Letting $\psi \colon G/\G\rightarrow \Tor^m$ denote the isomorphism between $G$ and $\Tor^m$, we have that $T_t$ is isomorphic to the flow $S_t=\psi T_t \psi^{-1}$ acting on $\Tor^m$ by translation by the one-parameter subgroup $\{\psi(a_t)\}$.

Write $\psi(a_t)=\textbf{b}_t = (b_{t,1}, \ldots, b_{t,m})\in \Tor^m$ for all $t\in \R$.  Then each $\brac{b_{t,i}}$ is a one-parameter subgroup of $\Tor$ and hence there is some $\a_i \in \R$ such that $b_{t,i}=\a_it $ for all $t\in\R$.  As $S_t$ is ergodic, $\brac{\a_1,\ldots,\a_m}$ are rationally independent, i.e., every non-trivial rational combination of $\a_1,\ldots,\a_m$ is non-zero.

It remains to show that for each $\bfx\in\Tor^m$,
\begin{eqnarray*} \{(S_{p_1(\bfs)}\bfx,\ldots,S_{p_k(\bfs)}\bfx)\}_{\bfs\in\R^d}& = & \{(x_1+p_1(\bfs)\alpha_1,\ldots,x_m+p_1(\bfs)\alpha_m,  \\
&  &\ldots , x_1+p_k(\bfs)\alpha_1, \ldots,x_m+p_k(\bfs)\alpha_m    )\}_{\bfs\in\R^d}  \end{eqnarray*}
is uniformly distributed in $\Tor^{km}$.  As the polynomials $\brac{\a_ip_j\colon 1\leq i \leq m  , 1\leq j \leq k }$ are rationally independent (i.e., every non-trivial rational combination of the polynomials $\brac{\a_ip_j\colon 1\leq i \leq m  , 1\leq j \leq k }$ is non-constant), this follows from Theorem \ref{btheorem5.2} below.\end{proof}

\begin{theorem}[Weyl, \cite{w}]\label{btheorem5.2}  Suppose $q_1,\ldots,q_w \colon \R^d \rightarrow \R$ are rationally independent polynomials.  Then $\{(q_1(\bfs),\ldots,q_w(\bfs))\}_{\bfs \in \R^d}$ is uniformly distributed in $\Tor^w$. 
\end{theorem}

We record the following consequence of the proof of Proposition \ref{bprop5.1} for future use:

\begin{proposition}\label{bprop5.3} Let $(X=G/\G,\sG/\G,\m)$ be a connected nilmanifold such that $G$ is abelian.  Then any nilflow on $X$ is isomorphic to translation by a one parameter subgroup on some finite dimensional torus.\end{proposition}

\begin{remark}\label{brmk5.1}  \emph{It is worth noting that Theorem \ref{btheorem1.2} fails if the polynomials $\{p_1,\ldots,p_k \}$ are not $\R$-independent.  Suppose there exist $a_1,\ldots, a_k \in \R$, not all zero, and $c\in \R$, so that $a_1p_1(\bfs) + \ldots + a_kp_k(\bfs)=c$ for all $\bfs\in\R^d$.  For each $i\in \{1,\ldots,k \}$, let $\{T_{a_i,t}\}_{t\in\R}$ be the flow on the torus $\Tor=\R/\Z$ defined by $T_{a_i,t}(x)=x+a_it$ for all $x\in \Tor$ and all $t\in\R$.  Let $S_t=T_{a_1,t} \times \ldots \times T_{a_k,t}$ and let $f_j(x_1,\ldots,x_k)=e^{2 \pi i x_j}\in L^{\infty}(\Tor^k)$ for all $j\in \{1,\ldots,k \}$.  Then  $$\avgRd S_{p_1(\bfs)}f_1 \cdot \ldots \cdot S_{p_k(\bfs)}f_k \,d\bfs$$ converges to $e^{2 \pi i (x_1 + \ldots + x_k+c)}$ in $L^2(\m)$ as $R_1,\ldots,R_d\rightarrow \infty$.}\end{remark}

\subsection{General description of the limit.}\label{bsec5.2}  In this section we compute the $L^2$-limit of \eqref{maineq}.  By \eqref{beq2.1} and Corollary \ref{bcor3.4} it suffices to compute \eqref{maineq} in the case where $(X=G/\G,\X,\m,\{T_t\})$ is a nilflow induced by some one-parameter subgroup $\{a_{t}\}_{t\in\R}$ of $G$.  We note that by Proposition \ref{bprop4.1}, in order to compute this limit, it suffices to describe for $x\in X$ the closure of the orbit
\begin{equation} \label{beq5.1} \{(a_{p_1(\bfs)}x,\ldots,a_{p_k(\bfs)}x)\}_{\bfs\in\R^d}\end{equation}
in $X^k$.  Leibman gives a description of orbits of the form \eqref{beq5.1} in \cite{Le4}.  In this section we show that in order to compute the limit of \eqref{maineq} it suffices to describe \eqref{beq5.1} when $p_1,\ldots,p_k$ are linear. 

\begin{proposition}\label{bprop5.4} Suppose $\{p_1,\ldots,p_k\colon \R^d\rightarrow \R\}$ is a collection of polynomials of the form $\{ \sum_{i=1}^{l}  \a_{1,i}q_i, \ldots, \sum_{i=1}^{l} \a_{k,i}q_i  \}$ for some collection of $\R$-independent polynomials $\{q_1,\ldots,q_l\colon \R^d \rightarrow \R\}$ with $q_i(0)=0$ for $i=1,\ldots,l$, and with $  \a_{j,i}\in\R$ for $i=1,\ldots,l$ and $j=1,\ldots,k$.  If $f_0,\ldots,f_k\in L^{\infty}(\m)$, then the averages 
\begin{equation*}\frac{1}{R^d} \int_0^R\cdots \int_0^R \int f_0 \cdot \prod_{j=1}^{k} T_{p_j(\bfs)}f_j  \,d\m \,d\bfs \end{equation*}
and 
\begin{equation*}\frac{1}{R^{l}}\int_{0}^{R} \ldots \int_{0}^{R} \int f_0 \cdot \prod_{j=1}^{k} T_{\sum_{i=1}^{l} \a_{j,i}u_i}f_{j} \,d\m\,d\bfu  \end{equation*}
have the same limit as $R \rightarrow \infty$.
\end{proposition}

A discrete time version of Proposition \ref{bprop5.4}, for averages along collections of three polynomials of Weyl complexity 2, is proved in \cite{Fran}.

\begin{proof}  We adapt the method of \cite{Fran} (Lemma 4.3).  By Corollary \ref{bcor3.4} and \eqref{beq2.1} it suffices to verify the lemma when the system is an ergodic nilflow, say $(X=G/\G,\sG/\G, \m, T_t)$, induced by some one-parameter subgroup $\{a_t\}$ of $G$.  By Proposition \ref{bprop4.4} it suffices to show that for every $x\in X$ the sets
  \begin{equation*}\label{eqn534}A=\{ (a_{u_0}x,a_{u_0+\sum_{i=1}^{l}  \a_{1,i}u_i}x, \ldots, a_{u_0+\sum_{i=1}^{l}  \a_{k,i}u_i}x)\}_{u_0,\ldots,u_l\in\R}  \end{equation*}
  and 
  \begin{equation*}\label{eqn535}B=\{ (a_{u_0}x,a_{u_0+\sum_{i=1}^{l}  \a_{1,i}q_i(\bfs)}x, \ldots, a_{u_0+\sum_{i=1}^{l}  \a_{k,i}q_i(\bfs)}x)\}_{u_0\in\R,\bfs\in \R^d}\end{equation*}
have the same closure.\footnote{The sets $A$ and $B$ are both subsets of $X^{k+1}$, despite the fact that $A$ is parameterized by $\R^{l+1}$ and $B$ is parameterized by $\R^{d+1}$.}  Identifying $X$ with $G^0/\G\cap G^0$, as in Section \ref{bsec4.2}, it suffices to assume $G$ is connected.
  
By Proposition \ref{bprop4.1} the closure of $A$ is a connected nilmanifold of the form $H /\Delta$, where $H$ is a connected closed subgroup of $G^{k+1}$ and $\Delta= H \cap \G^{k+1}$.  $B$ is clearly contained in $H/\Delta$ and it remains to be shown that $\overline{B}=H/\Delta$.  

Let $\pi \colon H/\Delta \rightarrow H/([H,H]\Delta)$ be the natural projection.  Then $\overline{\pi(A)}=H/([H,H]\Delta)$ and hence by Proposition \ref{bprop4.4} it suffices to show that $\overline{\pi(B)}=\overline{\pi(A)}$.  As $H$ is connected, Proposition \ref{bprop5.3} applies.  Thus we have reduced to showing that if $X=\Tor^m$, $\bm{\gamma} \in \Tor^m$, and the rotation $\bfx \mapsto \bfx+t \bm{\gamma} $ is ergodic, then for all $\R$-independent polynomials $q_1,\ldots,q_l\colon \R^d \rightarrow \R$, and for all $\bfx\in X$, the sets
\begin{equation}\label{beq5.2}\{ \bigl(  \bfx +u_0\bm{\gamma},\bfx+(u_0+\sum_{i=1}^{l}  \a_{1,i} u_i)\bm{\gamma}, \ldots, \bfx+( u_0+\sum_{i=1}^{l}  \a_{k,i} u_i)\bm{\gamma}\bigr)\}_{u_0,\ldots,u_l\in\R}  \end{equation}
  and 
  \begin{equation}\label{beq5.3}\{ \bigl(\bfx+u_{0}\bm{\gamma},\bfx+(u_0+\sum_{i=1}^{l}  \a_{1,i}q_i(\bfs))\bm{\gamma}, \ldots,\bfx+( u_0+\sum_{i=1}^{l}  \a_{k,i}q_i(\bfs))\bm{\gamma}\bigr)\}_{u_0\in\R,\bfs\in\R^d}\end{equation}
have the same closure.  

Write $\bm{\gamma}=(\gamma_1,\ldots,\gamma_m)$.  Let $\Lambda= \{1,\ldots,k\}\times\{1 ,\ldots, l\}\times\{1,\ldots,m \}$ and write $\widetilde{\a}_{\lambda}= \a_{j,i}\gamma_v$, $q_{\lambda}=q_i$, and $u_{\lambda}=u_i$ for all $\lambda=(j,i,v)\in\Lambda$.  Let $\Lambda'\subseteq \Lambda$ such that $\{\widetilde{\a}_{\lambda}q_{\lambda}(\bfs)\}_{\lambda\in\Lambda'}$ is maximal rationally independent subset of $\{\widetilde{\a}_{\lambda}q_{\lambda}(\bfs)\}_{\lambda\in\Lambda}$.  Let $w= \abs{\Lambda '}$ and write $\Lambda '=\{\lambda_1,\ldots,\lambda_w\}$.  Every element of $\{\widetilde{\a}_{\lambda}q_{\lambda}(\bfs)\}_{\lambda\in\Lambda}$ can be written as a linear combination of $\{\widetilde{\a}_{\lambda}q_{\lambda}(\bfs)\}_{\lambda\in\Lambda'}$ with rational coefficients:
\begin{equation*} \widetilde{\a}_{\lambda} q_{\lambda}(\bfs)=   c_{\lambda,\lambda_1}\widetilde{\a}_{\lambda_1}q_{\lambda_1}(\bfs)+ \ldots +  c_{\lambda,\lambda_m}\widetilde{\a}_{\lambda_w}q_{\lambda_w}(\bfs).
\end{equation*}
Write each $c_{\lambda,\lambda'}$ as a quotient of integers and let $d$ be the least common multiple of the denominators.  For each $\lambda \in \Lambda'$, let $\b_{\lambda}=\frac{\widetilde{\a}_{\lambda}}{d}$.  Then $\{\b_{\lambda}q_{\lambda}(\bfs)\}_{\lambda\in\Lambda'}$ is rationally independent and every element of $\{\widetilde{\a}_{\lambda}q_{\lambda}(\bfs)\}_{\lambda\in\Lambda}$ can be written as a linear combination of $\{\b_{\lambda}q_{\lambda}(\bfs)\}_{\lambda\in\Lambda'}$ with integer coefficients.

As the collections $\{q_1,\ldots,q_l\}$ and $\{u_1,\ldots,u_l\}$ are both $\R$-independent, they have the same dependence relations.  In particular, $\{\b_{\lambda}u_{\lambda}\}_{\lambda\in \Lambda'}$ is rationally independent.  By Theorem \ref{btheorem5.2}, $\{\bigl(\b_{\lambda_1}q_{\lambda_1}(\bfs) ,\ldots,  \b_{\lambda_w}q_{\lambda_w}(\bfs)\bigr) \}_{\bfs\in\R^d}$ and $\{\bigl(\b_{\lambda_1}u_{\lambda_1}, \ldots, \b_{\lambda_w}u_{\lambda_w}\bigr)\}_{u_{\lambda_1},\ldots,u_{\lambda_w}\in\R}$ are each equidistributed in $\Tor^w$.  Thus for each set of fixed values $u_{1},\ldots,u_{l}\in\R$, there exists $\bfs\in \R^d$ such that the distance in $\Tor^w$ between $\b_{\lambda}q_{\lambda}(\bfs) $ and $\b_{\lambda}u_{\lambda}$ is as small as desired for all $\lambda\in\Lambda'$.  If for every $\lambda\in\Lambda$
\begin{equation*} \widetilde{\a}_{\lambda}q_{\lambda}(\bfs)=   m_{\lambda,\lambda_1}\b_{\lambda_1}q_{\lambda_1}(\bfs)+ \ldots +  m_{\lambda,\lambda_w}\b_{\lambda_w}q_{\lambda_w}(\bfs)
\end{equation*}
for integers $\{m_{\lambda,\lambda_i}\}$, then for every $\lambda\in\Lambda$
\begin{equation*} \widetilde{\a}_{\lambda} u_{\lambda}=   m_{\lambda,\lambda_1}\b_{\lambda_1}u_{\lambda_1}+ \ldots +  m_{\lambda,\lambda_w}\b_{\lambda_w}u_{\lambda_w}.
\end{equation*}
Thus $\bfs\in\R^d$ can be chosen so that the distance in $\Tor^w$ between $\sum_{i=1}^{l}  \a_{j,i}q_i(\bfs)\bm{\gamma}$ and $\sum_{i=1}^{l}  \a_{j,i} u_i\bm{\gamma}$ is as small as desired for all $j\in\{1,\ldots,k\}$.  This shows that \ref{beq5.2} is contained in the closure of \ref{beq5.3}.  A similar argument shows that \ref{beq5.3} is contained in the closure of \ref{beq5.2}.
\end{proof}

By Proposition \ref{bprop5.4}, in order to describe the limit of \eqref{maineq} in general it remains to give a description of the limit along linear polynomials.  Let $G/\G$ be an $r$-step nilmanifold.  For each $k \in \{1,\ldots,r\}$, let $\Omega_k = \{(l_1,\ldots,l_m) \in \N^k \colon \sum_{i=1}^m l_i = k   \}$.  Given $\a_{j,i}\in\R$, for all $i\in\{1,\ldots,d\}$ and $j\in\{1,\ldots,k\}$, define the set

\begin{eqnarray*} H= \big\{\left( b_0 \prod_{n=1}^{r} \,\,\,\prod_{\om \in \Omega_{n}} b_{n,\om}^{\prod_{i=1}^{d}{\alpha_{1,i} \choose l_i}}, \ldots,  b_0 \prod_{n=1}^{r} \,\,\,\prod_{\om \in \Omega_{n}} b_{n,\om}^{\prod_{i=1}^{d}{\alpha_{k,i} \choose l_i}}\right) \colon \hspace{0.5in}  \\
\hspace{1.0in} b_{n,\om} \in G_n \,\forall \,n\in \{1,\ldots,r\},\forall \,\om \in \Omega_{n} \big\} \end{eqnarray*}
and let $\Delta= \G^k \cap H$.  $H$ is a closed subgroup of $G^k$, and the discrete subgroup $\Delta$ is cocompact \cite{Le4}.  Thus $H/\Delta$ is a nilmanifold with a Haar measure $m_H$.  

\begin{theorem}[Leibman, \cite{Le4}]\label{btheorem5.5}  Let $(X=G/\G, \sG/\G,\m, T_t)$ be an ergodic nilflow and let $\{p_1,\ldots,p_k \colon \R^d \rightarrow \R\}$ be a nice family of polynomials of the form $\{\sum_{i=1}^{d}\a_{1,i}s_i , \ldots,  \sum_{i=1}^{d}\a_{k,i}s_i\}$.  If $f_1,\ldots,f_k\in L^{\infty}(\m)$ then for a.e. $x=g\G \in X$ 
\begin{eqnarray}\nonumber \lefteqn{ \lim_{R_1,\ldots,R_d\rightarrow \infty} \avgRd f_1(T_{p_1(\bfs)}x) \cdot \ldots \cdot f_k(T_{p_k(\bfs)}x) \,d\bfs } \\ 
\nonumber & \hspace{1.5in}= & \int_{H/\Delta} f_1(gy_1\G) \cdot \ldots \cdot f_k(gy_k \G) \,dm_H(y\Delta),\end{eqnarray}
where $y=(y_1,\ldots,y_k)$, and $H$, $\Delta$ are as above.
\end{theorem}
Theorem \ref{btheorem5.5} follows from a specific case of Theorem 8.3 in \cite{Le4}, and the fact that for each $n\in\N$, the polynomials $\{\prod_{i=1}^{d}{s_i \choose l_i} \colon \om \in \Omega_{n}\}$ algebraically generate the polynomials $\R^d \rightarrow \R$ of degree $n$.  For further explaination, see \cite{Le4}, Section 0.7.  The discrete time version of Theorem \ref{btheorem5.5} in the $d=1$ case was given by Ziegler in \cite{Z2}.

\bigskip

\emph{\textbf{Limit Formula.}}  We now compute the $L^2$-limit of \eqref{maineq}.  If necessary, rewrite \eqref{maineq} so that $p_i(0)=0$ for $i=1,\ldots,k$.  Write $\{p_1,\ldots,p_k\}$ in the form $\{ \sum_{i=1}^{l}  \a_{1,i}q_i, \ldots, \sum_{i=1}^{l} \a_{k,i}q_i \}$, where $\{q_1,\ldots,q_l  \} $ are $\R$-independent polynomials.  By Theorem \ref{mainconvg}, the $L^2$-limit and the weak limit of \eqref{maineq} both exist and coincide.  Thus, by Proposition \ref{bprop5.4} the limit will be unchanged if we replace $\{ \sum_{i=1}^{l}  \a_{1,i}q_i, \ldots, \sum_{i=1}^{l} \a_{k,i}q_i \}$ with the linear polynomials $\{ \sum_{i=1}^{l}  \a_{1,i}u_i, \ldots, \sum_{i=1}^{l} \a_{k,i}u_i \}$.  Let $r\in \N$ such that $\sZ_r$ is characteristic for \eqref{maineq}.  After replacing $f_1,\ldots,f_k$ with their projections on $\sZ_r$ we assume that $\X=\sZ_r$.  As $\sZ_r$ is an inverse limit of r-step nilsystems, we can further assume that our system is an ergodic nilflow and compute the limit using Theorem \ref{btheorem5.5} (or by the more general method given by Theorem 8.3 in \cite{Le4}).

If $\sZ_1$ is characteristic then we can assume that our system is an ergodic flow given by multiplication by a one-parameter subgroup on a compact abelian Lie group $G$ with the Haar measure $\m$.  Identifying $X$ with $G^0/(\G \cap G^0)$, as in Section \ref{bsec4.2}, we may assume $G$ is connected, so $X=\Tor^m$ for some nonnegative integer $m$.  Thus for every $f_1,\ldots,f_k \in L^{\infty}(\m)$ the $L^2$-limit of \eqref{maineq} is 
\begin{equation}\label{beq5.4}\int_{\Tor^m}\ldots\int_{\Tor^m}  \prod_{j=1}^{k} f_{j}(x+\sum_{i=1}^{l} \a_{j,i}u_i) \,d\m(\bfu)\end{equation}
for a.e. $x\in\Tor^m$.

\subsection{Complexity.}\label{bsec5.3}  We define the \textbf{\emph{flow average complexity}} of a given family of polynomials $\{p_1,\ldots,p_k\}$ to be the smallest value of $r\in \N$ such that for any flow $(X,\X,\m,\{T_t\})$, $\sZ_r(\{T_t\})$ is characteristic for \eqref{maineq}.  We just write \textbf{\emph{complexity}} when it is clear we are referring to the flow average complexity.  For applications, it is useful to know the complexity of specific collections of polynomials.  By Proposition \ref{bprop5.4}, it suffices to compute the complexity for linear polynomials.  Combining Proposition \ref{bprop5.4} with Corollary \ref{bcor3.4}, we have the following result.  

\begin{corollary}\label{bcor5.6}  The complexity of a family $\{p_1,\ldots,p_k\}$ of non-constant essentially distinct polynomials is at most $k-1$.
\end{corollary}

A similar bound holds for the discrete average complexity in the case where the polynomials are all linear \cite{HK}.  However, it is still unknown whether in the discrete time setting a version of Corollary \ref{bcor5.6} holds for general families of polynomials.

It is shown in \cite{Furst1,CL,FWa,HK,HK2,L1} that if $\{p_1,\ldots,p_k \colon \Z^d \rightarrow \Z\}$ is a nice family of polynomials with $p_i(0)=0$ for all $i \in \{1,\ldots,k\}$, then there is some $r\in \N$ so that for each probability space $(X,\X,\m)$, for each measure preserving transformation $T \colon X \rightarrow X$, and for each F\o lner sequence $\folN$ in $\Z^d$, $\sZ_{r}(X,T)$ is \emph{characteristic} for the discrete time average
\begin{equation}\label{beq5.5}   \frac{1}{\abs{\Phi_N}} \sum_{\bfn \in \Phi_N} T^{p_1(\bfn)}f_1 \cdot \ldots \cdot T^{p_k(\bfn)}f_k,  \,\,\,\,f_1,\ldots,f_k \in L^{\infty}(\m)   .\end{equation}
In other words, for any $f_1,\ldots,f_k \in L^{\infty}(\m)$ with $\E(f_i | \sZ_{k-1})=0$ for some $i \in \{1,\ldots,k\}$, the average \eqref{beq5.5} converges to zero in $L^2(\m)$ as $N \rightarrow \infty$.  In this paper, we will refer to the minimal such $r\in\N$ as the \textbf{\emph{discrete average complexity}} of $\{p_1,\ldots,p_k \colon \Z^d \rightarrow \Z\}$.  A method for calculating the discrete average complexity is given in Section 6 of \cite{Le4}.

\begin{proposition}\label{bprop5.7}  Let $\{p_1,\ldots,p_k \colon \R^d \rightarrow \R\}$ be a family of linear polynomials with $p_i(0)=0$ and $p_i(\Z^d) \subseteq \Z$ for all $i\in\{1,\ldots,k\}$.  Then the flow average complexity of $\{p_1,\ldots,p_k \colon \R^d \rightarrow \R\}$ is bounded by the discrete average complexity of $\{p_1,\ldots,p_k \colon \Z^d \rightarrow \Z\}$.
\end{proposition}

\begin{proof}  Let $r$ be the discrete average complexity of $\{p_1,\ldots,p_k \colon \Z^d \rightarrow \Z\}$.  It suffices to show that if $f_1,\ldots,f_k\in L^{\infty}(X)$, with $\E(f_i|\sZ_r)=0$ for some $i \in \{1,\ldots,k\}$, then the $L^2$-limit of \eqref{maineq} is zero.

First suppose $T_1$ is totally ergodic.  If $\E(f_i|\sZ_r)=0$, then $\E(T_{p_i(\bfs)}f|\sZ_r)=0$ for all $\bfs\in\R^d$.  By an argument similar to the proof of Lemma \ref{blemma3.2}, and the Dominated Convergence Theorem,
 \begin{eqnarray}\label{beq5.6}\lefteqn{ \lim_{R_1,\ldots,R_d \rightarrow \infty} \frac{1}{\prod_{i=1}^{d}R_i} \int_{0}^{R_1}\ldots \int_{0}^{R_d}\prod_{i=1}^{k}   T_{p_i( \bfs)}f_i \, d\bfs }\\ 
\nonumber & = & \int_{[0,1]^d}  \lim_{R_1,\ldots,R_d \rightarrow \infty}\frac{1}{\prod_{i=1}^{d}\floor{R_i}}\sum_{n_1=0}^{\floor{R_1}-1} \ldots \sum_{n_d=0}^{\floor{R_d}-1}\prod_{i=1}^{k} T_1^{p_i( \bfn)}(T_{p_i( \bfs)}f_i) \, d\bfs . 
\end{eqnarray}

As $r$ is the discrete average complexity of $\{p_1,\ldots,p_k \colon \Z^d \rightarrow \Z\}$, the integrand is zero, and hence \eqref{beq5.6} is equal to zero.

If $T_1$ is not totally ergodic, fix $u \in \R$ such that $u>0$ and $T_u$ is totally ergodic.  Let $\{\widetilde{T}_t\}_{t\in\R}$ be the flow given by $\widetilde{T}_t= T_{ut}$ for all $t\in\R$.  Then $\widetilde{T}_{1}$ is totally ergodic and hence the $L^2$-limit of \eqref{maineq} is zero when $T_t$ is replaced with $\widetilde{T}_t$.  The change of variable $(s_1,\ldots,s_d) \mapsto (us_1,\ldots,us_d)$ now gives the result.  
\end{proof}

Combining Theorem \ref{btheorem1.2} and Remark \ref{brmk5.1}, we can characterize all families of complexity 0:
\begin{corollary}\label{bcor5.8}  A family $\{p_1,\ldots,p_k\}$ of non-constant essentially distinct polynomials has complexity 0 if and only if $\{p_1,\ldots,p_k\}$ are $\R$-independent.
\end{corollary}

\subsection{Bounding the complexity in examples.}\label{bsec5.4}  Let $\{p_1,\ldots,p_k\colon \R^l \rightarrow \R \}$ be a nice family of polynomials.  Define the \emph{\textbf{$p_j$-complexity}} of $\{p_1,\ldots,p_k\colon \R^l \rightarrow \R \}$ to be the smallest value of $r\in \N$ such that whenever $f_j\in L^{\infty}(\m)$ with $\E(f_j | \sZ_r)=0$, the average \eqref{maineq} converges to zero in $L^2(\m)$.  Then the complexity of $\{p_1,\ldots,p_k\colon \R^l \rightarrow \R \}$ is the maximum of the $p_j$-complexities for $j=\{1,\ldots,k\}$.  

We describe a method for determining a bound for the $p_1$-complexity.  Bounds for the other complexities can be determined by a similar process.  By Proposition \ref{bprop5.4}, it suffices to assume $\{p_1,\ldots,p_k\colon \R^l \rightarrow \R \}$ are linear polynomials of the form $\{ \sum_{i=1}^{l}  \a_{1,i}u_i, \ldots, \sum_{i=1}^{l} \a_{k,i}u_i \}$.  By relabeling the variables, we can further assume $\a_{1,1}\neq 0$.  Let $\Lambda_1=\{\a_{j,1}\colon \a_{j,1} \neq 0, 1 \leq j \leq k\}$.  Then $\Lambda_1$ are all coefficients of the variable $u_1$ in $\{ \sum_{i=1}^{l}  \a_{1,i}u_i, \ldots, \sum_{i=1}^{l} \a_{k,i}u_i \}$.  In Proposition \ref{bprop5.9} below, we show that if $\a_{1,1}\neq\a_{j,1}$ for all $2 \leq j \leq k$, then the $p_1$-complexity is at most $\abs{\Lambda_1}-1$.  Later in this section we explain how any collection of polynomials can be replaced by a collection of polynomials with the same complexities and which meets these requirements.  

\begin{proposition}\label{bprop5.9}  Suppose $\{p_1,\ldots,p_k\colon \R^l\rightarrow \R\}$ is a collection of distinct linear polynomials of the form $\{ \sum_{i=1}^{l}  \a_{1,i}u_i, \ldots, \sum_{i=1}^{l} \a_{k,i}u_i  \}$ with $  \a_{j,i}\in\R$ for $i=1,\ldots,l$ and $j=1,\ldots,k$, and let $\Lambda_1=\{\a_{j,1}\colon \a_{j,1} \neq 0, 1 \leq j \leq k\}$.  If $\a_{1,1}\neq 0$, and $\a_{1,1}\neq\a_{j,1}$ for all $2 \leq j \leq k$, then the $p_1$-complexity of $\{p_1,\ldots,p_k\colon \R^l\rightarrow \R\}$ is no greater than $\abs{\Lambda_1}-1$.
\end{proposition}

A similar type of result, for discrete time averages along collections of three polynomials of Weyl complexity 2, is proved in \cite{Fran}. 

\begin{proof}  We adapt the method of \cite{Fran} (Lemma 4.2).  Let $r=\abs{\Lambda_1}-1$, and for all $\bfu\in \R^l$ write $\bfu=(u_1,\ldots,u_l)$.  Let $f_1,\ldots, f_k\in L^{\infty}(\m)$ with $\norm{f_i}_{\infty} \leq 1$ for $i=1,\ldots,k$.  It suffices to show that if $\E(f_1|\sZ_r)=0$ then the $L^2$-limit of \eqref{maineq} is zero.  By Theorem \ref{mainconvg} the $L^2$-limit of \eqref{maineq} is identical to the $L^2$-limit of 
\begin{equation}\label{beq5.7}
\lim_{N \rightarrow \infty} \frac{1}{a(N) \cdot m(\bfRN)} \int_{\bfRN}\int_0^{a(N)} T_{p_1(\bfu)}f_1 \cdot \ldots\cdot T_{p_k(\bfu)}f_k \,d\bfu
\end{equation}
where $\bfRN=[-N,N]^{l-1}$ for all $N\in\N$ and $a(N)$ is an increasing sequence of integers to be chosen as follows.  By Corollary \ref{bcor3.4}, the $(\a_{j,1}u_1)$-complexity of the family $\{\a_{j,1}u_1\}_{j\in\Lambda}$ is at most $r$.  Write $\tilde{p}_j(u_2,\ldots,u_l)=\sum_{i=2}^{l} \a_{j,i}u_i$ for all $j=1,\ldots,k$ and note that if $\E(f_1|\sZ_r)=0$ then $\E(f_1 \circ T_{\tilde{p}_1(u_2,\ldots,u_l)}|\sZ_r)=0$ for all $u_2,\ldots,u_l \in \R$.  Since the map $\bfRN \rightarrow L^2(\m)$ given by $\tilde{\bfu}=(u_2,\ldots,u_l) \mapsto  \prod_{j \in \Lambda} f_{j}\circ T_{\tilde{p}_{j}(\tilde{\bfu})} $ is uniformly continuous, for each $N\in \N$ we are able to choose $a(N) \in \N$ with $a(N) > a(N-1)$ so that for all $\tilde{\bfu}=(u_2,\ldots,u_l)\in\bfRN$,
\begin{equation}\label{beq5.8}
\norm{\frac{1}{a(N)} \int_0^{a(N)} \prod_{j \in \Lambda}T_{\a_{j,1}u_1} \bigl(T_{\tilde{p}_{j}(\tilde{\bfu})}f_{j}\bigr)  \, du_{1}  }_{L^2(\m)} \leq \frac{1}{N}.
\end{equation}

Then for each $N\in \N$, 
\begin{align*}\label{eqn33}& \norm{\frac{1}{a(N) \cdot m(\bfRN)} \int_{\bfRN}\int_0^{a(N)} T_{p_1(\bfu)}f_1 \cdot \ldots\cdot T_{p_k(\bfu)}f_k \,d\bfu   }_{L^2(\m)}\\
&\leq  \frac{1}{m(\bfRN)} \int_{\bfRN}\norm{\frac{1}{a(N)} \int_0^{a(N)}  \prod_{j \in \Lambda} T_{\a_{j,1}u_1} \bigl( T_{\tilde{p}_{j}(\tilde{\bfu})}f_{j} \bigr)  \,du_1  }_{L^2(\m)} \,d\tilde{\bfu}.\end{align*}
By \eqref{beq5.8}, the $L^2$-limit of \eqref{beq5.7} is zero, which completes the proof.
\end{proof}

Now we describe how when $k\geq 2$, a nice family of linear polynomials $\{p_1,\ldots,p_k\colon \R^l \rightarrow \R \}$ can always be replaced with a collection of polynomials with the same complexities, which satisfies the requirements of Proposition \ref{bprop5.9}.

Define the coefficient matrix of $\{ \sum_{i=1}^{l}  \a_{1,i}u_i, \ldots, \sum_{i=1}^{l} \a_{k,i}u_i \}$ to be: 
\begin{equation*}A=\begin{pmatrix} 
\a_{1,1} & \a_{1,2} & \ldots \a_{1,l} \\
\a_{2,1} & \a_{2,2} & \ldots \a_{2,l}  \\
\vdots &  \vdots & \ldots & \vdots \\
\a_{k,1} & \gamma_{k,2} & \ldots \a_{k,l}.
\end{pmatrix}\end{equation*}
Notice that $\a_{1,1}$ will be distinct from $\{\a_{2,1},\ldots,\a_{k,1}\}$ if and only if the first entry of the first column is distinct from the other entries in that column.  By Proposition \ref{bprop5.4}, the $p_1$-complexity will remain unchanged if $u_1,\ldots,u_l$ are each replaced with $\R$-independent linear polynomials $q_1,\ldots,q_l$.  Write $q_j=\sum_{i=1}^l\gamma_{j,i}u_i$ for each $j\in \{1,\ldots,l\}$.  It is elementary to show that $q_1,\ldots,q_l$ are $\R$-independent if and only if the coefficient matrix $B=(q_{j,i})$ is invertible.  Moreover, the coefficient matrix $C$ of the produced polynomials will be the product, $C  =A \cdot B$.  Notice that the first entry of the first column of $C$ will be distinct from the other entries in that column if and only if $$\a_{j,1}\gamma_{1,1} + \a_{j,2}\gamma_{2,1} + \ldots + \a_{j,l}\gamma_{l,1}  \neq  \a_{1,1}\gamma_{1,1} + \a_{1,2}\gamma_{2,1} + \ldots + \a_{1,l}\gamma_{l,1} $$
 for all $j \in \{1,\ldots, k\}$ with $j\neq j'$.  This happens precisely when $(p_j-p_{1})(\bm{\gamma}) \neq 0$ for all $j \in \{2,\ldots, k\}$, where $\bm{\gamma} = (\gamma_{1,1},\ldots,\gamma_{l,1}) $.  As the solution set to the equation $(p_j-p_{1})(\bfx) = 0$ has measure zero for all $j\in \{2,\ldots, k\}$, there will certainly exist some non-zero $\bm{\gamma}$ with this property.  The remaining columns of $B$ can always be chosen so that $B$ is invertible.  Thus $B$ can always be found so that the resulting collection of polynomials will satisfy the hypotheses of Proposition \ref{bprop5.9}.  In many cases $B$ can be chosen so that the size of $\Lambda_1$ will be preserved, although it is unknown whether this will always be the case.  
 
\begin{example}\label{bex5.1}\emph{The collections of polynomials $\{u_1,2u_1,u_2\}$, $\{u_1,u_2,2u_1-u_2\}$, $\{u_1,u_2,u_3,\pi u_1+ \pi^2 u_3, 3u_2\}$, and $\{u_1, u_2,u_3, 2u_1+u_4,2u_2+u_4, 2u_3+u_4\}$ each have complexity at most 1, by a direct application of Proposition \ref{bprop5.9}.  None of these families are $\R$-independent, so by Corollary \ref{bcor5.8}, the complexity of each is 1. }
 \end{example}
 
\begin{example}\label{bex5.2}\emph{  The collection of polynomials $\{p_1(\bfu),p_2(\bfu),p_3(\bfu) \}=\{u_1,u_2,$ $u_1+u_2\}$ has complexity 1.  To see this, we use the method described above, setting $u_1=s+t$ and $u_2=s-t$.  The resulting collection of polynomials is $\{s+t,s-t,2s  \}$.  Using Proposition \ref{bprop5.9}, and examining the variable $t$, the $p_1$-complexity and the $p_2$-complexity are each at most 1.  Examining the variable $s$, we see that the $p_3$ complexity is also at most 1.  Taking the maximum of the $p_i$-complexities, we see that the complexity of $\{p_1(\bfu),p_2(\bfu),p_3(\bfu) \}$ is at most 1.  By Corollary \ref{bcor5.8}, the complexity is equal to 1.}
\end{example}

\begin{example}\label{bex5.3}\emph{  The collection of polynomials $\{p_1(\bfu),p_2(\bfu),p_3(\bfu),p_4(\bfu)  \}=$ $\{u_1,u_2,u_3, u_1+u_2+u_3\}$ has complexity 1.  To see this, use Proposition \ref{bprop5.4} and the change of variable $u_1=s$, $u_2=t$, $u_3=w-s$, to obtain the collection $\{s,t,w-s,t+w\}$ with the same complexities.  By proposition \ref{bprop5.9}, the $p_1$-complexity and the $p_3$-complexity are no greater than 1.  By symmetry, a similar change of variable shows the $p_2$-complexity is no greater than 1.  A different change of variable, $u_1=s$, $u_2=s+t$, and $u_3=w$, gives the collection $\{s,s+t,w,2s+t+w\}$, and by examining the coefficients of the variable $s$, we see that the $p_4$-complexity is at most 1.  By Corollary \ref{bcor5.8}, the complexity is 1.
}\end{example}

\begin{example}\label{bex5.4}\emph{  We show $\{p_1(\bfu),p_2(\bfu),p_3(\bfu),p_4(\bfu),p_5(\bfu),p_6(\bfu),p_7(\bfu)\}    = \{u_1,u_2,u_2,u_1+u_2,u_2+u_3,u_1+u_3\}$ has complexity at most 2.  By the change of variable $u_1=s$, $u_2=t-s$, $u_3=w-s$, we obtain the collection $\{s,t-s,w-s,t,w,t+w-2s,t+w-s\}$.  By examining the coefficients of the variable $s$, we see that the $p_1$ and $p_6$ complexities are at most 2.  By symmetry, the $p_2$, $p_3$, $p_4$, and $p_5$-complexities are at most 2.  The change of variable $u_1=s$, $u_2=t+s$, $u_3=w+s$, we obtain the collection $\{s,t+s,w+s, t+2s,w+2s,t+w+2s,t+w+3s\}$, and hence the $p_7$-complexity is bounded by 2.}
\end{example}

\begin{example}\label{bex5.5}\emph{  Let $l\geq 1$, and let $V = \{0,1\}^l$.  Let $P_l$ be the $l$-dimensional cube, i.e., $\{\bm{\varepsilon} \cdot \bfu\colon \bm{\varepsilon}\in V\}$.  The the complexity of $P_l$ is at most $l$.  This can be seen by doing a series of change of variables of the form $u_1 \mapsto u_1$, $(u_2,\ldots,u_l) \mapsto (u_2,\ldots,u_l) - u_1 \bm{\varepsilon}$, for $\bm{\varepsilon} \in \{-1,1\}^{l-1}$.}

\emph{Alternately, this fact follows from Proposition \ref{bprop5.7} and results in \cite{HK}.}
\end{example}

\begin{example}\label{bex5.6}\emph{  The collection $\{p_1(\bfu),p_2(\bfu),p_3(\bfu),p_4(\bfu)\}=\{\pi u_1+\pi^2u_2,$ $\pi^2u_1+\pi^3u_3,\pi u_1+\pi^2u_2+\pi u_3,\pi u_2+\pi u_3\}$ has complexity 1.  By the change of variable $(u_1,u_2,u_3) \mapsto (u_1,u_2+\frac{1}{\pi} u_1,u_3- \frac{1}{\pi}u_1)$, to get the collection $\{2\pi u_1+\pi^2  u_2, \pi^3 u_3, (2\pi-1)u_1+\pi^2 u_2+ \pi u_3,\pi u_2+\pi u_3 \}$, and examining the coefficients of the variable $u_1$, we see that the $p_1$ and $p_3$-complexities are also at most 1.  By the change of variable $(u_1,u_2,u_3) \mapsto (-\pi u_1,u_2+ u_1,u_3)$, to get the collection $\{\pi^2u_2,-\pi^3u_1+\pi^3u_3,\pi^2u_2+\pi u_3,\pi u_1+ \pi u_2+\pi u_3\}$, and examining the coefficients of the variable $u_1$, we see that the $p_2$ and $p_4$-complexities are also at most 1.  Thus the complexity is at most 1.  By Corollary \ref{bcor5.8}, the complexity is exactly 1.  }\end{example}

\begin{example}\label{bex5.7}\emph{The collection $\{t,2t,t^2\}$ has discrete average complexity 2 \cite{Fran}.  However, it is easily seen that the flow average complexity is 1.  To see this, by Proposition \ref{bprop5.4}, it suffices to examine the linear polynomials $\{u_1,2u_1,u_2\}$.  Thus by Proposition \ref{bprop5.9} the flow average complexity is at most 1, and is in fact equal to 1 by Corollary \ref{bcor5.8}. }\end{example}

It is unknown whether the bounds produced by the above method will achieve the flow average complexity for each nice family of polynomials.  

\section{Lower bounds.}\label{bsec6}  We now prove Theorems \ref{btheorem1.3} and \ref{btheorem1.4} using the method given in \cite{Fran}.
\begin{proof}[Proof of Theorem \ref{btheorem1.3}]  As much of this proof is identical to the proof of Theorem C (case 1) in \cite{Fran}, we give only a summary here.  

Without loss of generality, we assume $\{p_1,\ldots,p_k \}$ are non-constant and essentially distinct.  If $\{p_1,\ldots,p_k \}$ has complexity $0$, the result follows from Theorem \ref{btheorem1.2}.  

Suppose $\{p_1,\ldots,p_k \}$ has complexity $1$, and rewrite $\{p_1,\ldots,p_k \}$ in the form $\{q_1,\ldots,q_l, \sum_{i=1}^{l}  \a_{1,i}q_i, \ldots,\sum_{i=1}^{l} \a_{k-l,i}q_i \}$, for $\R$-independent polynomials $\{q_1,\ldots,q_l\}$.  By Proposition \ref{bprop5.3}, we may assume the Kronecker factor $\sZ_1$ is of the form $(\Tor^w, m, \{R_t\})$, where $R_t$ is defined by $R_t(\bfx)=\bfx+t\bm{\gamma}$ for some fixed $\bm{\gamma}\in \Tor^w$, and for all $\bfx \in \Tor^w$ and $t\in \R$.  Let $\pi_1 \colon X \rightarrow \sZ_1$ be the factor map.  For $\delta>0$, define the sets $V_{\delta}\colon= B(0,\delta)^l\subseteq \sZ_1^l$ and $S_{\delta}\colon = \{ \bfs\in \R^d \colon (q_1(\bfs)\bm{\gamma},q_2(\bfs)\bm{\gamma}, \ldots, q_l(\bfs)\bm{\gamma}) \in V_{\delta} \} $.  

First notice that $\sZ_1$ is characteristic for the average
\begin{equation}\label{beq6.1} \avgRd  1_{S_{\delta}}(\bfs) \cdot T_{p_1(\bfs)}f_1   \cdot \ldots \cdot T_{p_k(\bfs)} f_k \,d\bfs,\,\,\,f_1,\ldots,f_k\in L^{\infty}(\m).
\end{equation}
To see this, let $\chi_1,\chi_2,\ldots,\chi_l$ be any characters of $G$ and suppose $\E(f_i|\sZ_1)=0$ for some $i=1,\ldots,k$.  Then $\E(\chi_i\circ \pi_1 \cdot f_i | \sZ_1)=\chi_i \circ \pi_1\cdot \E(f_i  | \sZ_1)=0$, and hence
\begin{equation*}\avgRd T_{p_1(\bfs)}(\chi_1\circ \pi_1 \cdot f_1)    \cdot \ldots \cdot  T_{p_k(\bfs)}(\chi_k \circ \pi_1\cdot f_k) \,d\bfs.\end{equation*}
converges to zero in $L^2(\m)$ as $R_1,\ldots,R_d \rightarrow \infty$.  Approximating $1_{S_{\delta}}(\bfs)=$ $1_{V_{\delta}}(q_1(\bfs)\bm{\gamma},$ $\ldots,q_l(\bfs)\bm{\gamma})$ by functions of the form $\chi_1(q_1(\bfs)\bm{\gamma}) \cdot \ldots \cdot \chi_l(q_l(\bfs)\bm{\gamma})$, we see that $\sZ_1$ is characteristic for \eqref{beq6.1}.    

By Theorem \ref{btheorem5.2}, the path $\{(q_1(\bfs)\bm{\gamma},q_2(\bfs)\bm{\gamma},\ldots,q_l(\bfs)\bm{\gamma})  \}_{\bfs\in\R^d}$ is uniformly distributed in $\Tor^{wl}$, and hence
\begin{equation}\label{beq6.2}
\lim_{N \rightarrow \infty}  \frac{m\bigl(S_{\delta} \cap [0,R_1] \times \ldots \times [0,R_d]\bigr)}{R_1\cdot \ldots \cdot R_d} = m(V_{\delta}).
\end{equation}

It now follows from \eqref{beq5.4} and \eqref{beq6.2} that if $f_0,\ldots,f_k \in L^{\infty}(\m)$ and $\tilde{f}_i= \E(f_i | \sZ_1)$ for $i=1,\ldots,k$, then for any increasing sequence of rectangles $\folN$ in $\R^d$, each containing zero, with $\bigcup_{N \in\N}\Phi_N=\R^d$, we have
\begin{eqnarray}\nonumber\lefteqn{\lim_{N \rightarrow \infty} \frac{1}{m(S_{\delta} \cap \Phi_N)} \int_{S_{\delta} \cap \Phi_N} \int f_0 \cdot\prod_{j=1}^{k} T_{p_j(\bfs)}f_j \, d\m\,d\bfs}\\
\label{beq6.3} &=& \frac{1}{m(V_{\delta})} \int_{V_{\delta}} \int_G  \tilde{f}_0 \cdot\prod_{j=1}^{l} \tilde{f}_j(\bfx+u_j) \cdot \prod_{j=1}^{k-l}\tilde{f}_{l+j}(\bfx+\sum_{i=1}^{l} \a_{j,i}u_i) \,d\bfx\,d\bfu .\end{eqnarray} 

The limit of expression \eqref{beq6.3} as $\delta$ approaches zero is $\int \tilde{f_0} \cdot\tilde{f_1}\cdot\ldots\cdot\tilde{f_k} \,dm$.  Thus if $\delta$ is small enough and $f_i=f=\bm{1}_A$ for $i=0,1,\ldots,k$, then the quantity in \eqref{beq6.3} is greater than
$$ \int ( \tilde{f}\,)^{k+1} \,dm-\varepsilon \geq \Bigl(\int \tilde{f} \,dm\Bigr)^{k+1} -\varepsilon =\m(A)^{k+1}-\varepsilon  .$$
Therefore, if $\{ p_1,\ldots,p_k \}$ has complexity $1$, then for every $\varepsilon >0$ there exists $\delta>0$ so that 
$$ \lim_{N\rightarrow \infty } \frac{1}{m(S_{\delta}\cap\Phi_N)} \int_{S_{\delta}\cap\Phi_N} \m(A \,\cap\, T_{-p_1(\bfs)}(A) \,\cap\ldots \,\cap\, T_{-p_k(\bfs)}(A)) \,d\bfs$$
$$\geq \m(A)^{k+1} - \varepsilon .$$
\end{proof}
It is worth noting that it is our ability to give an explicit description of the limit of \eqref{maineq} in general which allows us to compute \eqref{beq6.1}, and hence to prove Theorem \ref{btheorem1.3} in its full generality. 

The proof of Theorem \ref{btheorem1.4} is identical to the proof of Theorem C (part 2) in \cite{Fran}, and thus we omit it.

\appendix \section{Appendix: The Correspondence Principle.}\label{appendixCP}

In this section we prove Theorem \ref{btheorem1.5} by modifying the proof of Proposition 2.2 in \cite{FKW}.

Let $E \subset \R$ such that $D^*(E)>0$.  Let $d\colon \R^2 \rightarrow \R$ be the Euclidean distance and define the function $\varphi \colon \R \rightarrow \R$ by
$$ \varphi(s)\colon = \min \{1, d(s,E)\}.$$

Let $X$ be the closure of the equicontinuous, uniformly bounded family of functions $\varphi_{t}(s)=\varphi(s+t)$ in the topology of uniform convergence over bounded sets in $\R$.  By the Ascoli-Arzel\'a Theorem, $X$ is compact.

We define a flow on $X$ by $T_{t}\psi(s)=\psi(s+t)$ for $\psi \in X$, $s,t\in\R$.  Since $D^*(E)>0$, there exists a sequence of intervals $S_n\subset \R$ such that
$$\frac{m(S_n \cap E)}{m(S_n)} \rightarrow D^*(E)>0,$$
and each interval $S_n$ induces a probability measure $\m_n$ on $X$:
$$\m_n (f)=\frac{1}{m(S_n)}\int_{S_n} f(T_{t}\varphi)\,dt.$$
By the Riesz Representation Theorem, Borel measures on $X$ correspond to linear functionals on $C(X)$, and thus there is a probability measure $\nu$ on $X$ and some subsequence $\{n_k\}$ such that
$$ \m_{n_k} \stackrel{\omega^*}{\rightarrow} \nu  .$$
Let $f_0 \colon X \rightarrow \R$ be the function given by $f_0(\psi)=\psi(0)$ for all $\psi \in X$.  Then $f_0$ is continuous.  Define $\widetilde{E}\subset X$ by 
$$ \psi \in \widetilde{E} \Leftrightarrow f_0(\psi)=0 \Leftrightarrow \psi(0)=0.   $$

\begin{lemma}\label{blemmaA.1} $\nu(\widetilde{E}) \geq D^*(E)$.
\end{lemma}
\begin{proof}  Recall that $\varphi(t)=0$ if $t\in E$.  For each $l\in\N$,
\begin{eqnarray*}\int_X (1-f_0(\psi))^l \,d\nu(\psi) & = &  \lim_{k \rightarrow \infty } \frac{1}{m(S_{n_k})}\int_{S_{n_k}} (1-f_0(T_t\varphi))^l \,dt\\
 & =  & \lim_{k \rightarrow \infty } \frac{1}{m(S_{n_k})}\int_{S_{n_k}} (1-\varphi(t))^l \,dt\\
& \geq & \lim_{k \rightarrow \infty}  \frac{m(S_{n_k}\cap E)}{m(S_{n_k})}  =D^*(E). \end{eqnarray*}
Thus $$\displaystyle{\nu(\widetilde{E} )= \lim_{l \rightarrow \infty } \int_X (1-f_0(\psi))^l \,d\nu(\psi)    \geq D^*(E)}.$$
\end{proof}
By the ergodic decomposition of $\nu$, there exists an ergodic measure $\m$ on $X$ such that $\m(\widetilde{E})\geq \nu(\widetilde{E})>D^*(E)$.  As $\sC(X)$ is separable, by the ergodic theorem $\m$-almost every $\psi\in X$ is a \textbf{generic point} for $\m$, i.e.,
$$ \lim_{R \rightarrow \infty} \frac{1}{R} \int_{0}^R f(T_{t}\psi)\,dt= \int f \,d\m $$
for every continuous function $f \in \sC(X)$.  

Furthermore, $\varphi \in X$ is \textbf{quasi-generic} (for the definition and proof of this fact in the discrete case, see \cite{F}) for $\m$, meaning there exists some sequence of intervals $\{I_N\}_{N \in \N}$ in $\R$ with ${\rm diam}(I_N) \rightarrow \infty$ such that
$$ \lim_{N \rightarrow \infty} \frac{1}{m(I_N)} \int_{I_N} f(T_t\varphi) \,dt =\int f \,d\m  $$
for every continuous function $f \in \sC(X)$.  To see that $\varphi$ is quasi-generic, let $\psi_0 $ be a generic point in $X$.  For each $f\in \sC(X)$,
$$ \lim_{N \rightarrow \infty}\frac{1}{N} \int_{0}^{N} f(T_t \psi_0) \,dt = \int f \, d\m    .$$
Let $\{f_k\}$ be a dense set of functions in $\sC(X)$, and let $N_k$ be an increasing sequence such that
$$ \abs{\frac{1}{N_k} \int_{0}^{N_k} f_j(T_t \psi_0) \,dt -\int f_j \, d\m           } <\frac{1}{2k} $$
for $j=1,2,\ldots,k$.  If $t_k$ is chosen so that the distance between $\psi_0$ and $T_{t_k}\varphi$ is sufficiently small, then
$$\abs{\frac{1}{N_k} \int_{t_k}^{N_k+t_k} f_j(T_{t} \varphi) \,dt -\int f_j \, d\m           }= \abs{\frac{1}{N_k} \int_{0}^{N_k} f_j(T_{t+t_k} \varphi) \,dt -\int f_j \, d\m           } <\frac{1}{k} $$
for $j=1,2,\ldots,k$.  Set $I_k=[t_k,N_k+t_k]$.  Then $\displaystyle{\lim_{N \rightarrow \infty} \frac{1}{m(I_N)} \int_{I_N} f_k(T_t\varphi) \,dt}$ $ =\int f_k \,d\m$ for each $k\in\N$.  By the density of $\{f_k\}$ in $\sC(X)$, $\varphi$ is quasi-generic for $\m$.

\begin{proposition}\label{bpropA.2}  Let $\m$ and $\tilde{E}$ be as above.  For $\{u_1,\ldots,u_l\}\subset \R$, we have for all $\delta >0$,
$$ D^*(\{ t\in\R\colon t, t+u_1,\ldots,t+u_l \in E_{\delta}\}) \geq \m(\tilde{E} \cap T_{u_1}^{-1} \tilde{E} \cap \ldots \cap T_{u_l}^{-1} \tilde{E})   .$$
\end{proposition}
\begin{proof}  Define the function $g\colon X \rightarrow \R$ by  
$$ g(\psi)= \begin{cases} \delta - f_0(\psi), & f_0(\psi)<\delta \\
 0, & f_0(\psi)\geq \delta \, . \end{cases} $$ 
As $\varphi$ is quasi-generic for $\mu$, there exists a sequence of intervals $I_N\subset\R$ such that 
$$\int f (\psi)\,d\m(\psi)  =  \lim_{N \rightarrow \infty}\frac{1}{m(I_N)} \int_{I_N}f(T_t\varphi) \,dt   $$
for all $f\in \sC(X)$.

Since $g(\psi)=\delta$ for $\psi \in \tilde{E}$, we have
\begin{eqnarray*}\lefteqn{ \delta^{l+1} \cdot \m(\tilde{E}\cap (T_{u_1}^{-1}\tilde{E}) \cap \ldots \cap (T_{u_l}^{-1}\tilde{E})) \leq  \int g(\psi) g(T_{u_1}\psi) \cdots g(T_{u_l}\psi) \,d\m(\psi)   }\\
& = & \lim_{N \rightarrow \infty}\frac{1}{m(I_N)} \int_{I_N} g(T_t \varphi) g(T_{u_1}T_t \varphi) \cdots g(T_{u_l}T_t \varphi) \,dt \hspace{0.9in}  \\
& \leq & \delta^{l+1} \cdot D^{*}( \{t\in\R \colon   t,t+u_1,\ldots,t+u_l \in E_{\delta}\})  .\end{eqnarray*}
\end{proof}

\section{ Appendix: van der Corput lemma.}\label{appendixVDC}

The following useful lemma is analogous to the discrete version given by van der Corput (see \cite{vdC}).
\begin{lemma}\label{vdc}  Let $(X,\m)$ be a probability space.  Suppose $(x,\bfs) \mapsto g_{\bfs}(x)$ is a map in $L^\infty(X \times \R^d)$ with $\norm{g_{\bfs}}_{L^\infty(\m)} \leq 1$ for almost every $  \bfs \in \R^d$.  Suppose $\nu$ is a Borel measure on $\R^d$ and let ${\Psi}$ be any $\nu$-measurable subset ${\Psi}\subseteq \R^d$ with $0<\nu({\Psi}) < \infty$.  Then
\begin{eqnarray*}\lefteqn{\label{vdceq}\limsup_{R_1,\ldots,R_d\rightarrow \infty} \norm{\avgRd g_{\bfs} \ d\bfs}_{L^2(\m)}^2  } \\
\nonumber & \leq  & \limsup_{R_1,\ldots,R_d\rightarrow \infty} \frac{1}{\nu({\Psi})^2} \int_{\Psi} \int_{\Psi} \avgRd \innerprod{g_{\bfs+\bfu},g_{\bfs+\bfv}}\, d\bfs \, d\bfu \, d\bfv . \end{eqnarray*}
\end{lemma}
\begin{proof}
Let ${\Psi} \subseteq \R^d$ with $0<\nu(\Psi)< \infty$, and let $R_N^{(1)},\ldots,R_N^{(d)}$ be sequences of positive real numbers with with $R_N^{(1)},\ldots,R_N^{(d)} \rightarrow \infty$ as $N \rightarrow \infty$.  Set $\Phi_N=[0,R_N^{(1)}] \times \ldots \times [0,R_N^{(d)}]$ for each $N\in \N$.  Then for all $N\in \N$, 
\begin{eqnarray*} \lefteqn{\avgNfol g_\bfs \ d\bfs  =  \frac{1}{\nu({\Psi})} \int_{\Psi} \avgNfol g_{\bfs}  \, d\bfs \, d\bfu} \\
& & = \!\frac{1}{\nu({\Psi})} \int_{\Psi} \frac{1}{m(\Phi_N)} \int_{\Phi_N }\! g_{\bfs+\bfu}  \, d\bfs \, d\bfu +\! \frac{1}{\nu({\Psi})} \int_{\Psi} \frac{1}{m(\Phi_N)} \int_{(\Phi_N -\bfu)\backslash \Phi_N }\!\!\! g_{\bfs+\bfu}  \, d\bfs \, d\bfu \\
& & \hspace{2.2in}  -\frac{1}{\nu({\Psi})} \int_{\Psi} \frac{1}{m(\Phi_N)} \int_{\Phi_N \backslash(\Phi_N -\bfu)}\! g_{\bfs+\bfu}  \, d\bfs \, d\bfu.\end{eqnarray*}
The last two terms approach zero as $N \rightarrow \infty$.  Thus, using the Cauchy-Schwarz Inequality, \eqref{vdceq} is equal to
\begin{eqnarray*} \lefteqn{  \limsup_{N\rightarrow \infty} \norm{ \frac{1}{\nu(\Psi)} \int_{\Psi} \avgNfol g_{\bfs+\bfu}  \, d\bfs \, d\bfu }_{L^2(\m)}^2 }\\
& \leq & \limsup_{N \rightarrow \infty} \frac{1}{m(\Phi_N)}\int_{\Phi_N} \norm{\frac{1}{\nu({\Psi})}\int_{\Psi}  g_{\bfs+\bfu} \, d\bfu}_{L^2(\m)}^2 \, d\bfs\\
& = & \limsup_{N \rightarrow \infty} \frac{1}{\nu({\Psi})^2} \int_{\Psi} \int_{\Psi} \avgNfol \innerprod{g_{\bfs+\bfu},g_{\bfs+\bfv}}\, d\bfs \, d\bfu \, d\bfv. \end{eqnarray*}
\end{proof}
We use the following corollaries of Lemma \ref{vdc}.
\begin{corollary}\label{bcorB.2}  Under the hypotheses of Lemma \ref{vdc},
\begin{eqnarray*}\lefteqn{\label{vdceq2}\limsup_{R_1,\ldots,R_d \rightarrow \infty} \norm{\avgRd g_{\bfs} \ d\bfs}_{L^2(\m)}^2  } \\
\nonumber & \leq  & \frac{1}{\nu({\Psi})^2} \int_{\Psi} \int_{\Psi}\limsup_{R_1,\ldots,R_d \rightarrow \infty} \abs{\avgRd \innerprod{g_{\bfs+\bfu},g_{\bfs+\bfv}}\, d\bfs} \, d\bfu \, d\bfv . \end{eqnarray*}
\end{corollary}
\begin{proof}  First use Lemma \ref{vdc}.  Then we are allowed to interchange the $\limsup$ and the integral by Fatou's Lemma.  \end{proof}
\begin{corollary}\label{bcorB.3}  Under the hypotheses of Lemma \ref{vdc}, there exists a sequence of rectangles $\Theta_N$ in $\R^{3d}$ with $\{0\} \subset \Theta_1 \subseteq \Theta_2\subseteq \Theta_3 \subseteq \ldots $, and $\bigcup_{N \in \N} \Theta_N=\R^{3d}$, such that 
\begin{eqnarray*} \lefteqn{ \limsup_{R_1,\ldots,R_d \rightarrow \infty}   \norm{\avgRd g_\bfs \ dm(\bfs)}_{L^2(\m)}^2 }  \\
& \hspace{0.5in}\leq  & \limsup_{N \rightarrow \infty} \frac{1}{m(\Theta_N)} \int_{(\bfs, \bfu, \bfv) \in\Theta_N}  \innerprod{g_{\bfs+\bfu},g_{\bfs+\bfv}} \, dm(\bfs) \, dm(\bfu) \, dm(\bfv).\end{eqnarray*}
\end{corollary}
\begin{proof}  
Let $\folN$ be an increasing sequence of rectangles $\folN$ in $\R^d$ with $\bigcup_{N\in\N}\Phi_N=\R^d$ and $0\in\Psi_N$ for all $N\in\N$, such that \begin{eqnarray}\nonumber\lefteqn{\limsup_{N \rightarrow \infty}   \norm{\avgNfol g_\bfs \, dm(\bfs)}_{L^2(\m)}^2 }\\
\label{J} & = & \limsup_{R_1,\ldots,R_d \rightarrow \infty}   \norm{\avgRd g_\bfs \, dm(\bfs)}_{L^2(\m)}^2   .\end{eqnarray}  
Let $J$ denote the quantity \eqref{J}.  Choose any increasing sequence of rectangles $\brac{\Psi_N}_{N \in \N}$ in $\R^d$ with $\bigcup_{N\in\N}\Psi_N=\R^d$ and $0\in\Psi_N$ for all $N\in\N$.  Using Lemma \ref{vdc}, find a sequence $\{ M_N\}_{N\in \N}\subseteq \N$ so that for each $N \in \N_+$, $M_N \geq N$ and
$$ \frac{1}{m(\Psi_N)^2} \int_{\Psi_N} \int_{\Psi_N} \frac{1}{m(\Phi_{M_N})} \int_{\Phi_{M_N}} \innerprod{g_{\bfs + \bfu},g_{\bfs+\bfv}} \, d\bfs \, d\bfu \, d\bfv > J-\frac{1}{N}.$$
Define $\Theta_N=\Phi_{M_N} \times \Psi_N \times \Psi_N$.  Then
$$ \limsup_{N \rightarrow \infty} \frac{1}{m(\Theta_N)} \int_{(\bfs,\bfu,\bfv) \in \Theta_N} \innerprod{g_{\bfs+\bfu},g_{\bfs+\bfv}} \, d\bfs \, d\bfu \, d\bfv            \geq J.$$
\end{proof}

\textbf{Acknowledgments.}  I would like to thank B. Kra for her guidance and N. Frantzikinakis and A. Leibman for valuable comments on the preprint.

\nocite{FK3,HK2,L3,Shah2,Bour}
\small

\bibliographystyle{amsplain} 
\bibliography{bibliographyA}
\end{document}